\documentclass[12pt]{amsart}
\usepackage{amssymb}
\usepackage{graphics}
\usepackage{colordvi}

\copyrightinfo{}{}
\newtheorem{thm}{Theorem}
\newtheorem{lemma}[thm]{Lemma}
\newtheorem{cor}[thm]{Corollary}
\newtheorem{prop}[thm]{Proposition}
\newtheorem{definition}[thm]{Definition}
\newtheorem{rem}[thm]{Remark}
\newtheorem{ex}[thm]{Example}

\newtheorem*{thm*}{Theorem}

\newenvironment{remark}{\begin{rem}\rm}{\end{rem}}

\newcommand{\ann}{\text{ann}}
\newcommand{\CC}{\mathbb{C}}

\newcommand{\OO}{\mathbb{O}}
\newcommand{\PP}{\mathbb{P}}
\newcommand{\RR}{\mathbb{R}}
\newcommand{\QQ}{\mathbb{Q}}
\newcommand{\OP}{\mathbb{O}\mathbb{P}}

\newcommand{\frb}{\mathfrak{b}}
\newcommand{\frg}{\mathfrak{g}}
\newcommand{\frh}{\mathfrak{h}}
\newcommand{\frl}{\mathfrak{l}}
\newcommand{\frk}{\mathfrak{k}}
\newcommand{\frn}{\mathfrak{n}}
\newcommand{\frp}{\mathfrak{p}}
\newcommand{\frq}{\mathfrak{q}}
\newcommand{\frr}{\mathfrak{r}}
\newcommand{\frs}{\mathfrak{s}}

\newcommand{\frz}{\mathfrak{z}}

\newcommand{\calO}{\mathcal{O}}

\newcommand{\Hom}{{\rm Hom}}

\newcommand{\into}{\hookrightarrow}
\newcommand{\onto}{\twoheadrightarrow}
\newcommand{\codim}{\mathop{\rm codim}}
\newcommand{\Gr}{{\rm Gr}}

\newcommand{\ten}{\mbox{\small $1\!0$}}
\newcommand{\eleven}{\mbox{\small $1\!1$}}

\newcommand{\Inv}{{\rm Inv}}
\newcommand{\Invc}{{\rm Inv}^c}
\newcommand{\Stab}{{\rm Stab}}
\newcommand{\stab}{{\rm stab}}
\newcommand{\opp}{o}

\newcommand{\pr}{{\rm pr}}

\begin{document}

\title[The recursive nature of cominuscule Schubert calculus]{The
  recursive nature of \\cominuscule Schubert calculus}  

\author{Kevin Purbhoo}
\address{Department of Combinatorics \& Optimization \\
         University of Waterloo \\
         Waterloo, ON, N2L 3G1\\
         CANADA}
\email{kpurbhoo@math.uwaterloo.ca}
\urladdr{http://www.math.uwaterloo.ca/\~{}kpurbhoo}
\author{Frank Sottile}
\address{Department of Mathematics\\
         Texas A\&M University\\
         College Station\\
         TX \ 77843\\
         USA}
\email{sottile@math.tamu.edu}
\urladdr{http://www.math.tamu.edu/\~{}sottile}
\thanks{Work of Sottile supported by the Clay Mathematical Institute,  NSF
         CAREER grant DMS-0538734, and Peter Gritzmann of the Technische
         Universit\"at M\"unchen} 
\thanks{Work of Purbhoo supported by NSERC} 

%
%
\begin{abstract}
 The necessary and sufficient Horn inequalities which determine 
 the non-vanishing 
 Littlewood-Richardson coefficients in the cohomology of a Grassmannian are
 recursive in that they are naturally indexed by non-vanishing
 Littlewood-Richardson coefficients on {\em smaller}\/ Grassmannians.
 We show how 
 non-vanishing in the Schubert calculus for cominuscule flag varieties
 is similarly recursive.
 For these varieties, the non-vanishing of products of Schubert classes is
 controlled by the non-vanishing products on smaller cominuscule flag
 varieties.
 In particular, we show that the lists
 of Schubert classes whose product is non-zero 
 naturally correspond to the integer points in 
 the feasibility polytope, which is defined by inequalities 
 coming from non-vanishing products of Schubert classes on smaller
 cominuscule flag varieties.
 While the Grassmannian is cominuscule, our necessary and sufficient
 inequalities are different than the classical Horn inequalities. 
\end{abstract}
\maketitle

\section*{Introduction}


We investigate the following general problem:
Given Schubert subvarieties $X, X',\dotsc, X''$ of a
flag variety, when is the intersection of their general
translates 
 \begin{equation}\label{E:generic}
   gX \cap g'X'\cap\dotsb\cap g''X''
 \end{equation}
non-empty?
When the flag variety is a Grassmannian, it is known that 
such an intersection is
non-empty if and only if the indices of the Schubert varieties,
expressed as partitions, satisfy the linear Horn inequalities.
The Horn inequalities are themselves 
indexed by lists of partitions
corresponding to such non-empty intersections on smaller
Grassmannians.
This recursive answer to our original question
is a consequence of work of Klyachko~\cite{Klyachko} who 
linked eigenvalues of sums of hermitian matrices, highest weight modules
of $\mathfrak{sl}_n$, and the Schubert calculus, and of
Knutson and Tao's proof~\cite{KT99} of 
Zelevinsky's Saturation Conjecture~\cite{Ze99}. 
These two results proved Horn's Conjecture~\cite{Ho62}
about the eigenvalues of sums of Hermitian matrices.
This had wide implications in mathematics (see the
surveys~\cite{Fu98,Fu00}) and 
raised many new and evocative questions. 
For example, the recursive nature of this geometric question
concerning the intersection of Schubert varieties was initially
mysterious, as the proofs used much
more than the geometry of the Grassmannian. 

Belkale~\cite{Be02} provided a geometric proof of the Horn
inequalities, which explains their recursive nature.
His method relied upon an analysis of the tangent spaces to Schubert varieties.
One of us (Purbhoo) reinterpreted Belkale's proof~\cite{Pu04} using
two-step partial flag varieties (Grassmannians are one-step partial
flag varieties) for the general linear group.
This approach starts from the observation that 
the non-emptiness of an intersection~\eqref{E:generic} can be translated into a question
of transversality involving the tangent spaces of Schubert varieties
(Proposition~\ref{P:gen-feasible}).

For other groups, two-step partial flag varieties are replaced by 
fibrations of flag varieties. 
Suppose that $R\subset P$ are parabolic subgroups of a complex
reductive algebraic group $G$.
Then $P/R=L/Q$, where $L$ is the Levi subgroup of $P$ and $Q$ is a parabolic subgroup of
$P$ and we have the fibration sequence of flag varieties.
 \[
   \begin{matrix}
   L/Q\ =\  P/R &\longrightarrow& G/R\\
        &&\raisebox{4pt}{\big\downarrow}\rule{0pt}{15pt}\\
        &&G/P
   \end{matrix}
 \]
Given Schubert varieties $X$ on $G/P$ and $Y$ on $L/Q$, there is a unique lifted Schubert
variety $Z$ on $G/R$ which maps to $X$ with fiber $Y$ over the generic point of $X$.
Each tangent space of $G/R$ has a map to $\frz$, the dual of the center of 
the nilradical of $R$.
Let $C(X,Y)$ be the codimension in $\frz$ of the image of the tangent 
space to $Z$ at a smooth point.

Suppose that we have Schubert varieties $X, X',\dotsc,X''$ of $G/P$ 
such that the
intersection~\eqref{E:generic} of their general translates is non-empty.
Given Schubert varieties $Y, Y',\dotsc,Y''$ of $L/Q$ whose general translates (by elements
of $L$) have non-empty intersections, then we have the inequality 
 \begin{equation}\label{Eq:simplified_ineqs}
   C(X,Y)+C(X',Y')+\dotsb+C(X'',Y'')\ \leq\ \dim \frz\,.
 \end{equation}
We show that a subset of these necessary inequalities are
sufficient to determine when a general intersection~\eqref{E:generic} of Schubert
varieties is non-empty, when $G/P$ is a \Blue{cominuscule flag variety}.
For each cominuscule $G/P$, we identify a set 
$M(P)$ of parabolic subgroups   
$Q\subset L$.   
We state a version of our main result (Theorem~\ref{Theorem}).

\begin{thm*}
 Suppose that $G/P$ is a cominuscule flag variety.
 Then the intersection~\eqref{E:generic} is non-empty if and only if 
 for every $Q\in M(P)$ and every Schubert varieties $Y,Y',\dotsc,Y''$ of $L/Q$
 whose general translates have non-empty intersection,
 the inequality~\eqref{Eq:simplified_ineqs} holds.
\end{thm*}

As discussed in Section~\ref{S:Two}, this solves the question of when
an arbitrary product of Schubert classes on a cominuscule flag variety 
is non-zero. 

The subgroups $Q\in M(P)$ have the property that $L/Q$ is also cominuscule,
and thus the inequalities which determine the non-emptiness of~\eqref{E:generic} are
recursive in that they come from similar non-empty intersections on smaller cominuscule
flag varieties. 
For Grassmannians, these inequalities are different than the Horn 
inequalities, and hence give a new proof of the Saturation Conjecture.
Moreover, the inequalities for the Lagrangian and orthogonal Grassmannians
are different, despite their having the same sets of solutions!

By cominuscule flag variety, we mean the orbit of a highest weight vector
in (the projective space of) a cominuscule representation of 
a linear algebraic group $G$.
These are analogs of the Grassmannian for other Lie types;
their Bruhat orders are distributive lattices~\cite{P84} and the 
multiplication in their cohomology rings is governed by a uniform 
Littlewood-Richardson rule~\cite{TY07}.
Cominuscule flag varieties $G/P$ are distinguished in that the unipotent
radical of $P$ is abelian~\cite{RRS92} and in that a Levi subgroup
$L$ of $P$ acts on the tangent space at $eP$ with finitely
many orbits.
There are other characterizations of cominuscule flag varieties
which we discuss in Section~\ref{S:cominuscule}.
We use that $G/P$ is 
cominuscule in many essential ways in our arguments,
which suggests
that cominuscule flag varieties are the natural largest
class of flag varieties for which these tangent space methods can be
used to study the non-vanishing of intersections~\eqref{E:generic}.

Since the algebraic groups 
$G$ and $L$ need not have the same Lie 
type, in many cases the
necessary and sufficient 
inequalities of Theorem~\ref{Theorem} are indexed by non-empty
intersections of Schubert varieties 
on cominuscule flag varieties of a 
different type.  
For example, the inequalities for the Lagrangian Grassmannian
are indexed by non-empty intersections on ordinary Grassmannians.  
This is in contrast to the classical Horn recursion, which is purely in 
type $A$, involving only ordinary Grassmannians.  
Thus the recursion we
obtain is a recursion within the class of cominuscule flag varieties,
rather than a type-by-type recursion.
This is reflected in our proof of the 
cominuscule recursion, which is entirely independent of type;
in particular we do not appeal to the
classification of cominuscule flag varieties.

This paper is structured as follows.
Section~\ref{S:One} establishes our notation and develops background material.
Section~\ref{S:Two} states our main theorem 
 precisely (Theorem~\ref{Theorem})
and derives necessary inequalities 
 (Theorem~\ref{T:necessary}), which are more general than
the inequalities~\eqref{Eq:simplified_ineqs}.
Section~\ref{S:three} contains the proof of our main theorem, some of which
relies upon technical results about root systems, which are given in the
Appendix.
In Section~\ref{S:four} we examine the cominuscule recursion  
in more detail, describing it on a case-by-case basis.
In Section~\ref{S:five} we compare our results and inequalities to other
systems of inequalities for non-vanishing in the Schubert calculus,
including the classical Horn inequalities, and the dimension inequalities
of Belkale and Kumar~\cite{BK04}.
We have attempted to keep Sections~\ref{S:three} and~\ref{S:four}
independent, so that the reader who is more interested in examples 
may read them in the opposite order.

\section{Definitions and other background material}\label{S:One}
We review basic definitions and elementary facts that we use concerning
linear algebraic groups, Schubert varieties and their tangent spaces, transversality, and
cominuscule flag varieties.
All algebraic varieties, groups, and algebras will be over the complex numbers,
as the proofs we give of the main results are valid only for complex varieties.

\subsection{Linear algebraic groups and their flag varieties}
We assume familiarity with the basic theory of algebraic groups and Lie
algebras as found in~\cite{Bo91,FH91,Hu72,Se66}.
We use capital letters $B,G,H,K,L,P,Q,R,\dotsc$ for algebraic
groups and the corresponding lower-case fraktur letters for their Lie
algebras $\frb,\frg,\frh,\frk,\frl,\frp,\frq,\frr,\dotsc$.
We also use lower-case fraktur letters $\frs,\frz$ for subquotients
of these Lie algebras.
Throughout, $G$ will be a reductive algebraic group, $P$ a parabolic
subgroup of $G$, $B\subset P$ a Borel subgroup of $G$, and $e\in G$ will be the
identity. 
Let $H$ be a maximal torus of $G$ with $H\subset B$.
Let $L\subset P$ be the Levi (maximal reductive) subgroup containing $H$.
We have the Levi decomposition $P=LN_P$ of $P$ with $N_P$ its 
($H$-stable) unipotent radical.
Write $G^{ss}$ and $L^{ss}$ for the semisimple parts of $G$ and $L$,
respectively.
Write $W$ or $W_G$ for the Weyl group of $G$, which is the quotient
$N_G(H)/H$. 
Note that $W_P=W_L$.

There is a dictionary between parabolic subgroups $Q$ of $L$ and parabolic
subgroups $R$ of $P$ which contain a maximal torus of $L$, 
\[
   Q\ =\ R\cap L\quad\mbox{ and }\quad R\ =\ Q N_P\,. 
\]
Thus $R$ is the maximal subgroup of $P$ whose restriction to $L$ is $Q$.
We will always use the symbols $Q$ and $R$ for parabolic subgroups of 
$L$ and $P$ associated in this way.
We will typically have $H\subset Q (\subset R)$.
Set $B_L:=B\cap L$, a Borel subgroup of $L$ that contains $H$.
We say that $Q$ and $R$ are \Blue{{\em standard parabolic subgroups}} 
if $B_L \subset Q$ (equivalently $B \subset R$).
Then the surjection $\pr\colon G/R\onto G/P$ has fiber $P/R= L/Q$, so we
have the fibration diagram:
 \begin{equation}\label{E:FlagFibration}
   \begin{matrix}
    L/Q\ =\ P/R &\longrightarrow& G/R\\
        &&\raisebox{9pt}{\Big\downarrow\makebox[0.01in][l]{\pr}}\rule{0pt}{22pt}\\
        &&G/P
   \end{matrix}
 \end{equation}

Let $\Phi\subset\frh^*$ be the roots of the Lie algebra $\frg$.
These decompose into positive and negative roots,
$\Phi=\Phi^+\sqcup\Phi^-$, where $\Phi^-$ are the roots of $\frb$.
Our convention that the roots of $B$ are negative will simplify the
statements of our results. 
Write $\Delta$ for the basis of simple roots in $\Phi^+$.
For $\alpha\in\Phi$, let $\frg_\alpha$ be the (1-dimensional) 
$\alpha$-weight space of $\frg$.
Then we have 
\[
   \frg\ =\ \frh\oplus\bigoplus_{\alpha\in\Phi}\frg_\alpha\qquad
   \mbox{and}\qquad
   \frb\ =\ \frh\oplus\bigoplus_{\alpha\in\Phi^-}\frg_\alpha\ .
\] 
We write $\Phi(\frs)$ for the non-zero weights of an $H$-invariant subquotient
$\frs$ of $\frg$, and 
$\Phi^+(\frs)$ for $\Phi(\frs)\cap\Phi^+$.  
 Note that these weights are all roots of $\frg$.
The Killing form on $\frg$ pairs $\frg_\alpha$ with $\frg_{-\alpha}$
and identifies $\frg$ with its dual. 
Under this identification, 
the dual $\frs^*$ of an $H$-invariant subquotient $\frs$ 
is another subquotient of $\frg$, and 
$\Phi(\frs^*) = - \Phi(\frs)$.  
In this way, the dual of $\frn_P$ is identified with $\frg/\frp$.

The Weyl group $W$ acts on all these structures.
For example, if $g\in N(H)$, then the conjugate $gBg^{-1}$ of
$B$ depends only upon the coset $gH$, 
which is 
the element $w$ of $W$ determined by $g$.
Write $wBw^{-1}$ for this conjugate, and 
use similar notation for conjugates of other subgroups of $G$.
Conjugation induces a left action on roots and we have $w\Phi^-=\Phi(wBw^{-1})$. 
The \Blue{{\em inversion set}} of $w\in W$ is the set of positive roots which
become negative under the action of $w$,
$\Inv(w):=w^{-1}\Phi^-\cap \Phi^+$.
The inversion set determines $w$, and the cardinality of $\Inv(w)$
is the \Blue{{\em length}} of $w$, $\ell(w):=|\Inv(w)|$.

Borel subgroups containing $H$ are conjugate by elements of $W$.
For 
$w\in W_G$, $wBw^{-1}\cap P$ is a solvable subgroup of
$P$ which is not necessarily maximal.
However, $w Bw^{-1}\cap L$ is a Borel subgroup of $L$, and this
has a nice description in terms of the Weyl groups $W_G$ and $W_L=W_P$.
Let $\pi\in w W_L$ be the  coset representative of minimal length (with
respect to reflections in the simple roots $\Delta$).
Write $W^P$ for this set of minimal length coset representatives, and similarly
write $W^Q$ for the set of minimal length representatives of cosets of $W_Q$ in
$W_L$. 
Set $\lambda:=\pi^{-1}w\in W_L$.
Then  $\ell(w)=\ell(\lambda)+\ell(\pi)$.
This corresponds to a decomposition of the inversion set of $w$.
Note that $\Phi^+ = \Phi^+(\frl) \sqcup \Phi(\frg/\frp)$.
Then
 \begin{equation} \label{E:inversionequation} 
 \begin{split}
  \Inv(\lambda)\ &=\  \Inv(w)\cap \Phi(\frl)\,,\\
  \Inv(\pi)    \ &=\ \lambda \Inv(w)\cap\Phi(\frg/\frp)\,,\quad\mbox{and}\\
  \Inv(w)      \ &=\  \Inv(\lambda)\sqcup \lambda^{-1} \Inv(\pi)\,.
 \end{split}
 \end{equation}

\subsection{Schubert varieties and their tangent spaces}
\label{S:schubertvarieties}
Points of the flag variety $G/P$ are parabolic subgroups conjugate to
$P$, with $gP$ corresponding to the subgroup $gPg^{-1}$.
A Borel subgroup $B$ of $G$ acts with finitely many orbits on $G/P$.
When $H\subset B\subset P$, each orbit has the form $BwP$ for some $w\in W$.
The coset $wP$ is the unique $H$-fixed point in the orbit $BwP$.

If $wW_P=w'W_P$ for some $w,w'\in W$, then $wP=w'P$.
Thus these $B$-orbits are indexed by the set $W^P$.  
If $P'\in B\pi P$ for $\pi\in W^P$, then we say that $P'$ has 
\Blue{{\em Schubert position}} $\pi$
with respect to the Borel subgroup $B$.
When this happens, there is a $b\in B$ such that $bP'b^{-1}\supset \pi B\pi^{-1}$.  
The decomposition 
\[
  G/P\ =\ \coprod_{\pi\in W^P} B \pi P
\]
of $G/P$ into $B$-orbits is the Bruhat decomposition of $G/P$.
The orbit $X^\circ_\pi B:= B\pi P$ 
is called a \Blue{{\em Schubert cell}}
and is parametrized by the unipotent subgroup
$B\cap \pi N_{P^\opp}\pi^{-1}$, where $N_{P^{\opp}}$ is
the unipotent radical of the parabolic subgroup $P^\opp$
opposite to $P$.
The closure of $X^\circ_\pi B$
is the \Blue{{\em Schubert variety}} $X_\pi B$, 
which has dimension $\ell(\pi)$.

For each $\pi\in W^P$, define the \Blue{{\em planted Schubert cell}}
$X^\circ_\pi$ to be the translated orbit $\pi^{-1} B\pi P$, and the 
\Blue{{\em planted Schubert variety}}
$X_\pi$ to be its closure.
A translate of the Schubert cell $X^\circ_\pi B$ contains $eP$ if and 
only if it 
has the form $p X^\circ_\pi$ for some $p\in P$.

The tangent space to $G/P$ at $eP$ is naturally identified with the
Lie algebra quotient $\frg/\frp$.  
As the nilpotent subgroup $N_{P^\opp}$
parameterizes $G/P$ in a neighborhood of $eP$, the
tangent space can also be 
identified with $\frn_{P^\opp}$;
indeed, the natural map,
$\frn_{P^\opp} \to \frg/\frp$, is an $H$-equivariant isomorphism.
Since $X^\circ_\pi$ is parametrized by 
$\pi^{-1}B\pi\cap N_{P^o}$, its tangent space \Blue{$T_\pi$} at $eP$ 
(an $H$-submodule of $\frn_{P^\opp}$) has weights
 \begin{align*}
  \Phi(T_\pi)\ 
  & =\ (\pi^{-1}\Phi^-)\cap \Phi(\frn_{P^o})  \\
  & =\ (\pi^{-1}\Phi^-)\cap \Phi(\frg/\frp) \ =\ \Inv(\pi)\,.
 \end{align*}
$P$ acts on the tangent space $T_{eP}G/P = \frg/\frp$.  Translating
$T_\pi \subset \frg/\frp$ by $p \in P$, we obtain the
tangent space $p T_\pi$ to $p X_\pi$ at $eP$.

These planted Schubert varieties and their tangent spaces fit into the fibration
diagram~\eqref{E:FlagFibration}.
Let $R\subset P$ be standard parabolic subgroups of $G$, 
and set $Q:=L\cap R$ be the standard parabolic
subgroup of $L$ corresponding to $R$.
Minimal coset representatives of $W_R$ in $W_G$ are products 
$\pi\lambda$, 
where $\pi\in W^P$ and $\lambda\in W^Q$ is a minimal representative of
$W_Q$ in $W_L$.
Then the image of the Schubert cell $B\pi\lambda R$ of $G/R$ in $G/P$ is the
Schubert cell $B\pi P$.
When $\pi$ is the identity, we have that 
$B \lambda R/R=B_L \lambda Q/Q$.

In general, the fiber $B\pi \lambda R\to B\pi P$ is isomorphic 
to $B_L\lambda Q$.  In particular, we have 
\[
   \begin{matrix}
    X^\circ_\lambda &\longrightarrow& X^\circ_{\pi \lambda}\\
        &&\raisebox{9pt}{\Big\downarrow\makebox[0.01in][l]{\pr}\rule{0pt}{16pt}}\\
        &&\lambda^{-1} X^\circ_\pi
   \end{matrix}
\]
and thus we obtain a short exact sequence of the tangent spaces
 \begin{equation}\label{Eq:Tan_SES}
    T_\lambda\ \hookrightarrow\ T_{\pi\lambda}\ \twoheadrightarrow\ 
    \lambda^{-1}T_\pi\,.
 \end{equation}
Indeed, if $bRb^{-1}\in X^\circ_{\pi\lambda}$ lies in the fiber, then $b$ lies in
$\lambda^{-1}\pi^{-1}B\pi\lambda \cap P$.
Since $R$ contains the unipotent radical of $P$,
we can assume that in fact
\[
  b\ \in\ \lambda^{-1}\pi^{-1} B \pi \lambda \cap L\ =\ \lambda^{-1} B_L \lambda
\]
and thus $bRb^{-1}\cap L=bQb^{-1}\in  X^\circ_\lambda$.
The converse is straightforward.
Here, we used that $\pi^{-1}B\pi \cap L=B_L$, which follows from
$\Inv(\pi)\cap\Phi(\frl)=\emptyset$.

\subsection{Transversality}
We write $V^*$ for the linear dual of a vector space $V$ and write
$U^\ann$ for the annihilator of a subspace $U$ of $V$.
A collection of linear subspaces of $V$
\Blue{{\em meets transversally}} if their annihilators are in direct sum.

For us, a variety will always mean a reduced, but not-necessarily
irreducible scheme over the complex numbers.
A collection of algebraic subvarieties of a smooth variety $X$ is
\Blue{{\em transverse at a point $p$}} if they are each smooth at
$p$ and if
their tangent spaces at $p$ meet transversally, as subspaces of the
tangent space of $X$ at $p$.
A collection of algebraic subvarieties of a smooth variety $X$ 
\Blue{{\em meets transversally}} if they are transverse at 
the generic point of every component in their intersection.
We freely invoke Kleiman's Transversality Theorem~\cite{Kl74}, which asserts
that if a (complex)
reductive algebraic group acts transitively on 
a smooth variety $X$, then general
translates of subvarieties of $X$ meet transversally.

We establish the following result from elementary linear algebra,
which will be indispensable in analyzing the transversality 
of Schubert varieties.


\begin{prop}\label{P:Transversality}
  Suppose that we have a short exact sequence of vector spaces
\[
   0 \ \longrightarrow\ W\ \longrightarrow\ V\ \longrightarrow\ V/W
   \ \longrightarrow\ 0\,.
\]
 Let $U_1,\dotsc,U_s$ be linear subspaces of\/ $V$ and set $S_i:=W\cap
 U_i$ and $M_i:=(S_i+W)/W$, for $i=1,\dotsc,s$.
\begin{enumerate}
 \item[(i)] If $U_1,\dotsc,U_s$ are transverse in $V$, then 
  $M_1,\dotsc,M_s$ are transverse in $V/W$.  
 \item[(ii)] If $S_1,\dotsc,S_s$ are transverse in $W$, then $U_1,\dotsc,U_s$
   are transverse if and only if $M_1,\dotsc,M_s$ are transverse.
\end{enumerate} 
\end{prop}

\begin{proof}
 It suffices to prove this for $s=2$, 
 as subspaces are transverse if and only if they are pairwise transverse.
 If $U^\ann_1,U^\ann_2$ form a direct sum, then their subspaces
 $M^\ann_1,M^\ann_2$ form a direct sum, and (i) follows immediately.
 This proves one direction of (ii).
 For the other, consider its dual statement:
 If $S^\ann_1+S^\ann_2$ and $M^\ann_1+M^\ann_2$ are
 direct sums, then so is $U^\ann_1+U^\ann_2$.
 Note that $M^\ann_i=U^\ann_i\cap(V/W)^*$ and 
 $S^\ann_i$ is the image of $U^\ann_i$ in $W^*$.
 But if $U^\ann_1+U^\ann_2$ is not a direct sum, then 
 $U^\ann_1\cap U^\ann_2\neq\{0\}$.
 By assumption on $M^\ann_1$ and $M^\ann_2$, the image of 
 $U^\ann_1\cap U^\ann_2$ in $W^*$ is a
 non-empty subspace lying in $S^\ann_1\cap S^\ann_2$.
\end{proof}

It follows immediately from the definition of transversality that
if $U_1, \dots, U_s$ are transverse linear subspaces of $V$, then we
must have the codimension inequality
\[
 \sum_{i=1}^s \codim U_i\ \leq\ \dim V\,.
\]
We freely make use of this basic fact in our arguments.

\subsection{Cominuscule flag varieties}\label{S:cominuscule}

We list several equivalent characterizations of 
cominuscule flag varieties $G/P$.
Recall that $P=LN_P$ is the Levi decomposition of $P$.
Then
 \begin{enumerate}
  \item[(i)] $N_P$ is abelian.
  \item[(ii)] $L$ has finitely many orbits on $N_P$, equivalently on its Lie algebra
         $\frn_P$ and on $\frg/\frp=T_{eP}G/P$.
  \item[(iii)] $\frg/\frp$ is an irreducible representation of $L$, which 
     implies that the Weyl group $W_L$ acts transitively on roots of the 
     same length 
     in $\Phi(\frg/\frp)$. 
  \item[(iv)] $P=P_\alpha$ is a maximal parabolic subgroup of $G$ and 
     the omitted simple root $\alpha$ occurs with coefficient 1 in the 
     highest root of $G$.
 \end{enumerate}
Sources for these equivalences, with references, may be found
in~\cite{LM03,RRS92,Roe93}. 
Cominuscule flag varieties come in five infinite families 
with two exceptional cominuscule flag varieties.

Let $G/P$ be a cominuscule flag variety and $\alpha$ the root corresponding to
the maximal parabolic subgroup $P$.
As explained in~\cite{LM03}, 
the semisimple part $L^{ss}$ of the Levi subgroup of $P$
has Dynkin diagram obtained from that of $G$ by deleting the node
corresponding to the root $\alpha$.
The representation of $L^{ss}$ on the tangent space $\frg/\frp$
is the tensor product of fundamental representations given by
marking the nodes in the diagram of 
$L^{ss}$ that were adjacent to $\alpha$.
This is summarized in Table~\ref{Ta:One}.
\begin{table}[htb]

\begin{tabular}{|c|c|c|c|c|}\hline
 $G/P$ & $\Gr(k,n+1)$ & $Q^{2n-1}$ & $LG(n)$ & $Q^{2n-2}$
   \rule{0pt}{12pt}\\\hline
 $ G$  & $A_n$      & $B_n$      &  $C_n$  &  $D_n$\rule{0pt}{12pt}\\\hline
\raisebox{5pt}{$\alpha$}
         &\includegraphics{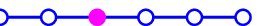}&
          \includegraphics{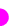}&
          \includegraphics{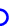}&
          \includegraphics{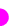}\rule{0pt}{19pt}\\\hline
$L^{ss}$    & $A_{k-1}\times A_{n-k}$ & $B_{n-1}$& $A_{n-1}$ &$D_{n-1}$%
  \raisebox{-3pt}{\rule{0pt}{15pt}}\\\hline
\raisebox{4pt}{$\frg/\frp$}
        &\includegraphics{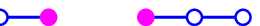}&
         \includegraphics{figures/TBn.eps}&
         \includegraphics{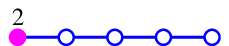}&
         \includegraphics{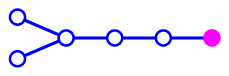}\\\hline
\end{tabular}\bigskip
\\
\begin{tabular}{|c|c|c|c|}\hline
 $G/P$ & $OG(n)$ & $\OP^2$ & $G_\omega(\OO^3,\OO^6)$\rule{0pt}{12pt} \\\hline
 $ G$  & $D_n$      & $E_6$      &  $E_7$\rule{0pt}{12pt}\\\hline
\raisebox{5pt}{$\alpha$}
       &  \includegraphics{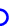}\ &
          \includegraphics{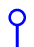}&
          \includegraphics{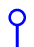}\rule{0pt}{19pt}\\\hline
$L^{ss}$    & $A_{n-1}$ & $D_5$& $E_6$\rule{0pt}{12pt} \\\hline
\raisebox{4pt}{$\frg/\frp$}
        &\includegraphics{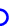}&
         \includegraphics{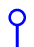}&
         \includegraphics{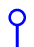}\\\hline
\end{tabular}\bigskip

\caption{Cominuscule Flag Varieties}
\label{Ta:One}
\end{table}

 The varieties $Q^{2n-1}$ and $Q^{2n-2}$ are odd- and even-dimensional 
 quadrics respectively.  $LG(n)$ is the Lagrangian Grassmannian.
The superscript 2 in the Dynkin diagram of 
$A_{n-1}$ in the column for $LG(n)$ indicates that this representation 
has highest weight twice the corresponding fundamental weight.
The second cominuscule flag variety in type $D_n$ is the \Blue{{\em orthogonal
    Grassmannian}}, $OG(n)$.
This is one of two components of the space of maximal isotropic 
subspaces in the vector space $\CC^{2n}$, which is equipped with a nondegenerate
symmetric bilinear form.
It is also known as the spinor variety.
The notation $\OP^2$ is for the Cayley plane (projective plane for the octonians)
and $G_\omega(\OO^3,\OO^6)$ is borrowed from~\cite{LM03} (as was the idea for
Table~\ref{Ta:One}). 

\section{Feasibility and Statement of Main Theorem}\label{S:Two}

The general problem that we are investigating is, given
 $\pi_1,\dotsc,\pi_s$ with $\pi_i \in W^P$
and general translates 
$g_1X_{\pi_1} B, \dotsc, g_sX_{\pi_s} B$ of the corresponding 
Schubert varieties, when is the intersection
 \begin{equation}\label{E:General_Intersection}
   g_1X_{\pi_1} B\ \cap\ g_2X_{\pi_2} B\ \cap\ \dotsb\ \cap\  g_sX_{\pi_s} B
 \end{equation}
non-empty?
A list $\pi_1,\dotsc,\pi_s$ with $\pi_i\in W^P$ is a 
\Blue{{\em Schubert position}} for $G/P$. 
It is \Blue{{\em feasible}} if 
such general intersections~\eqref{E:General_Intersection} are non-empty.
For $g\in G$, the translate $g X_\pi B$ is another Schubert variety,
but for the Borel subgroup $gBg^{-1}$.
Thus, $\pi_1,\dotsc,\pi_s$ is a feasible Schubert position
if, for any 
Borel subgroups $B_1,\dotsc,B_s$, there is a parabolic
subgroup $P'$ having Schubert position $\pi_i$ with respect to $B_i$ for
each $i=1,\dotsc,s$.

Feasibility is often expressed in terms of algebra.
Write $\sigma_{\pi}$ for the class of a Schubert variety $X_\pi B$ in the
cohomology ring of $G/P$.  Then the product 
$\prod_{i=1}^s \sigma_{\pi_i}$
is non-zero if and only if a general
intersection of the form~\eqref{E:General_Intersection} is non-empty,
if and only if the Schubert
position $\pi_1,\dotsc,\pi_s$ is feasible.
If $\sum_{i=1}^s \codim X_{\pi_i}B =  \dim G/P$, then 
the generic intersection~\eqref{E:General_Intersection} is finite,
and the integral
 \[
    \int_{G/P} \sigma_{\pi_1}\sigma_{\pi_2} \dotsb \sigma_{\pi_s}
 \]
computes the number of points in this intersection.
In this case we say that $\pi_1, \dotsc, \pi_s$ is a \Blue{{\em top-degree}}
Schubert position.

In this section, we state two theorems, Theorem~\ref{T:necessary} and our
main result, Theorem~\ref{Theorem}, which give conditions for 
feasibility in terms of inequalities.  We then show how the problem
of feasibility can be reformulated in terms of transversality 
for tangent spaces to Schubert varieties.  Using this, we prove
Theorem~\ref{T:necessary}.  The ideas in this section form the 
foundation for the proof of Theorem~\ref{Theorem}, which is given in 
Section~\ref{S:three}.

\subsection{Statement of main results}

 As in Section~\ref{S:schubertvarieties},
 let $R \subset P$ be standard parabolic subgroups of $G$, and 
 let $Q:=L \cap R$.
 Let $\frs$ be any $R$-submodule of 
 the nilradical $\frn_R$ of $\frr$.
 As $\frn_R^*$ is identified with the tangent space to $G/R$ 
 at $eR$, 
 dual to the inclusion $\frs\hookrightarrow\frn_R$ is the 
 surjection
\[
  \varphi_\frs\ \colon\ T_{eR}G/R\ \relbar\joinrel\twoheadrightarrow\ 
   \frs^*\,.
\]

\begin{thm}\label{T:necessary}
 Suppose that $\pi_1,\dotsc,\pi_s$ is a feasible Schubert position for $G/P$.
 Given any feasible Schubert position $\lambda_1,\dotsc,\lambda_s$ for $L/Q$,
 we have the inequality
 \begin{equation}
 \label{E:generalineq}
   \sum_{i=1}^s \codim \varphi_\frs(T_{\pi_i\lambda_i})\ 
     \leq\ \dim \frs\ .
 \end{equation}
\end{thm}

We prove Theorem~\ref{T:necessary} in Section~\ref{S:derivation}.

\begin{remark}
\label{R:combinatorial}
Note that each 
inequality~\eqref{E:generalineq} 
is a combinatorial condition:
As $T_{\pi_i\lambda_i}$ is $H$-invariant,
the left hand side can be calculated explicitly using
\[
  \codim \varphi_\frs(T_{\pi_i\lambda_i})\ =\  
  |\Phi(\frs^*) - \Phi(T_{\pi_i\lambda_i})| \ =\ 
  |\Phi(\frs^*) - \Inv(\pi_i\lambda_i)|\,.
\]
\end{remark}

As $Q$, $\frs$, and $\lambda_1,\dotsc,\lambda_s$ range over all possibilities,
this gives a system of necessary inequalities for the feasible 
Schubert position $\pi_1,\dotsc,\pi_s$.

The inequalities of Theorem~\ref{T:necessary} are more general than
those given in the introduction.  They specialize to a number of previously
known inequalities, which we discuss further in Section~\ref{S:five}.
For our main theorem, we
specialize to the case where $\frs = Z(\frn_R)$,
the center of the nilradical of $\frr$. 
In this case we write $\varphi_R$ for $\varphi_{Z(\frn_R)}$.  
Then 
the inequality~\eqref{E:generalineq} becomes
 \begin{equation}\label{E:w-ineq}
    \sum_{i=1}^s \codim\varphi_R(T_{\pi_i\lambda_i})\ 
    \leq\ \dim Z(\frn_{R})\,.
 \end{equation}

 Suppose $G/P$ is a cominuscule flag variety.
 Let $M(P)$ be the set of standard parabolic subgroups 
 of $L$ which are equal to the stabilizer of
 the tangent space (at some point) to some $L$-orbit 
 on $\frg/\frp$.
 We will show (Lemma~\ref{L:Qcominuscule}) that if $Q\in M(P)$, then $L/Q$ is
 cominuscule; however, not all parabolic subgroups 
 $Q$ of $L$ with $L/Q$ cominuscule
 are conjugate to a subgroup in $M(P)$
 (see Sections~\ref{S:An} and~\ref{S:Dn-spinor}).

We now state our main theorem, which is proved in Section~\ref{S:three}.

\begin{thm}\label{Theorem}
 Suppose that the semisimple part of $G$ is simple 
 (see Remark~\ref{R:Gsimple}).  
 Let $\pi_1,\dotsc,\pi_s$ be a Schubert position for a cominuscule flag
 variety $G/P$.
 Then  $\pi_1,\dotsc,\pi_s$ is feasible if and only if the following
 condition holds:
 for every $Q \in M(P) \cup \{L\}$ and 
     every feasible top-degree Schubert position
     $\lambda_1,\dotsc,\lambda_s$ for $L/Q$, 
     the inequality~\eqref{E:w-ineq} holds.
\end{thm}

The degenerate case of $Q=L$ in~\eqref{E:w-ineq} gives the
\Blue{{\em basic codimension inequality}}
 \begin{equation}\label{E:degreeineq}
   \sum \codim T_{\pi_i}\ \leq\ \dim G/P \,.
 \end{equation}
If we restrict our attention
to top-degree Schubert positions $\pi_1,\dotsc,\pi_s$, this degenerate 
case is unneeded
as~\eqref{E:degreeineq} is then an equality.  
Thus we have the following
recursion purely for the feasible top-degree Schubert positions.

\begin{cor}\label{C:topdegree}
 Suppose that the semisimple part of $G$ is simple.
 Let $\pi_1,\dotsc,\pi_s$ be a top-degree Schubert position for a 
 cominuscule flag variety $G/P$.
 Then $\pi_1,\dotsc,\pi_s$ is feasible if and only if for 
 every $Q\in M(P)$ and  every feasible top-degree Schubert position 
 $\lambda_1,\dotsc,\lambda_s$ for $L/Q$,
 the inequality~\eqref{E:w-ineq} holds.
\end{cor}

\begin{remark}
\label{R:Gsimple}
The hypothesis that $G^{ss}$ be simple is technically necessary, but
mild.
Theorem~\ref{Theorem} and Corollary~\ref{C:topdegree} allow us to obtain
necessary and sufficient
inequalities for any reductive group $G$ and cominuscule $G/P$.  
When $G^{ss}$ is not simple, the group $G/Z(G)$
is the product 
$G^1 \times \dots \times G^k$ of simple groups, and $P/Z(G)$
is the product $P^1 \times \dots \times P^k$ of parabolic subgroups
$P^j \subset G^j$.
Then $G/P \cong G^1/P^1 \times \dots \times G^k/P^k$, where each 
$G^j/P^j$ is cominuscule (or $P^j = G^j$). 
Furthermore, each Schubert position
$\pi_i \in W^P$ is a $k$-tuple 
$(\pi_i^1, \dotsc, \pi_i^k) \in W^{P^1} \times \dots \times W^{P^k}$, 
and $\pi_1, \dotsc, \pi_s$ is feasible for $G/P$
if and only if $\pi_1^j, \dotsc, \pi_s^j$ is feasible for $G^j/P^j$
for all $j$.  
Thus we simply check that
each $\pi_1^j, \dots, \pi_s^j$ satisfies the 
inequalities~\eqref{E:w-ineq} with $Q \in M(P^j) \cup \{L^j\}$ for all $j$.  

These inequalities are not of the form~\eqref{E:w-ineq} on $G/P$, but
rather of the more general form~\eqref{E:generalineq} on $G/P$.
\end{remark}

\begin{remark}\label{Rem:saturation}
In \cite{Be02}, Belkale showed that the Horn recursion
implies Zelevinsky's Saturation Conjecture.  As we will see in 
Sections~\ref{S:An} and~\ref{S:Horn}, our recursion for Grassmannians
is different from the classical Horn recursion.  Nevertheless, Belkale's 
argument can be used to show that 
our recursion also implies the Saturation Conjecture.
We will not repeat the argument here, but the reader who is familiar
with it will see that little modification is required.
Thus our proof of Theorem~\ref{Theorem} will implicitly also give a new 
proof of the Saturation Conjecture.  
\end{remark}

\begin{remark}
 As can be seen from the examples in Sections~\ref{S:even-quadric}
 and~\ref{S:odd-quadric}, the system of inequalities in 
 Theorem~\ref{Theorem} may be redundant.  An interesting problem is 
 to find an irredundant subset of these inequalities which solves the 
 feasibility question.
 For the classical Horn inequalities, this is known \cite{Be01, KTW04}, 
 however since our inequalities are different, this problem is open, 
 even for the Grassmannian.
\end{remark}

\subsection{Local criteria for feasibility}

The derivation of necessary inequalities of Theorem~\ref{T:necessary}
begins with the observation that feasibility can be detected locally.
Recall that $P$ acts on the tangent space $T_{eP} G/P \simeq \frg/\frp$.

\begin{prop}
\label{P:gen-feasible}
A Schubert position $\pi_1,\dotsc,\pi_s$ for $G/P$ is feasible 
if and only if the intersection 
 \begin{equation}\label{E:tanInt}
   p_1T_{\pi_1}\ \cap\  p_2T_{\pi_2}\ \cap\ \dotsb\ \cap\  p_sT_{\pi_s}
 \end{equation}
is transverse in $\frg/\frp$, for general $p_1,\dotsc,p_s\in P$.
\end{prop}

\begin{proof}
Since a general intersection
of Schubert varieties is transverse 
at the generic point of each of its components, either 
a general intersection 
 is empty or else it is (i) non-empty, (ii) of the expected dimension,
and (iii) the Schubert varieties meet transversally at every 
such generic point.
Thus, given an intersection~\eqref{E:General_Intersection}
which is non-empty but otherwise general, {\em either} it is transverse at
the generic point of every component
and the Schubert position is feasible, 
{\em or else} it is not transverse at  
the generic point of some component
and the Schubert position is infeasible.

Consider an intersection of
Schubert varieties~\eqref{E:General_Intersection} that are general subject
to their containing the distinguished point $eP$.  
 Such an intersection is of the form
 \begin{equation}
 \label{E:intplanted}
    p_1X_{\pi_1}\ \cap\  p_2X_{\pi_2}\ \cap\ \dotsb\ \cap\  p_sX_{\pi_s}\,,
 \end{equation}
 where $p_1, \dotsc, p_s$ are general elements of $P$.
Since $G/P$ is a homogeneous space, 
a general intersection~\eqref{E:intplanted} containing $eP$
is transverse
if and only if a non-empty but otherwise general
intersection~\eqref{E:General_Intersection} is transverse.
 But~\eqref{E:tanInt} is 
 just the intersection of the tangent spaces at $eP$ to the 
 Schubert varieties in~\eqref{E:intplanted}.
Thus the intersection~\eqref{E:tanInt} is transverse if and
only if $\pi_1, \dotsc, \pi_s$ is feasible.
\end{proof}

When $G/P$ is a cominuscule flag variety, we have 
the following refinement of 
Proposition~\ref{P:gen-feasible}, in which the
general elements $p_1,\dotsc,p_s\in P$ are replaced by general
elements $l_1,\dotsc,l_s\in L$ in~\eqref{E:tanInt}.

\begin{prop}\label{Prop:feasible}
 A Schubert position $\pi_1,\dotsc,\pi_s$ for a cominuscule flag variety $G/P$ is
 feasible if and only if the intersection 
\[
    l_1T_{\pi_1}\ \cap\  l_2T_{\pi_2}\ \cap\ \dotsb\ \cap\  l_sT_{\pi_s}
\]
 is transverse for generic $l_i\in L$.
\end{prop}

\begin{proof}
Since $G/P$ is cominuscule, the unipotent radical $N_P$ of $P$ is
abelian and thus acts  trivially on its Lie algebra $\frn_{P}$ and on
its dual, $\frg/\frp$. 
Thus we may replace general elements $p_1,\dotsc,p_s\in P$ by general
elements $l_1,\dotsc,l_s\in L$ in~\eqref{E:tanInt}.
\end{proof}

\subsection{Derivation of necessary inequalities}
\label{S:derivation}

\begin{prop}\label{P:FeasibleSP}
 For each $i=1,\dotsc,s$, let $\pi_i$ and $\lambda_i$ be Schubert positions for
 $G/P$ and $L/Q$, respectively, and $\pi_i\lambda_i$ the corresponding Schubert
 position for $G/R$. 
 \begin{enumerate}
  \item[(i)] If $\pi_1 \lambda_1,\dotsc,\pi_s \lambda_s$ is feasible, then 
      so is $\pi_1,\dotsc,\pi_s$.
  \item[(ii)] If both $\lambda_1,\dotsc,\lambda_s$ and
     $\pi_1,\dotsc,\pi_s$ is feasible, 
     then $\pi_1 \lambda_1,\dotsc,\pi_s \lambda_s$ is feasible.
 \end{enumerate}
\end{prop}

\begin{proof}
 For (i), the hypotheses imply that on $G/R$ the intersection
 \begin{equation}\label{E:genint_G/R}
   g_1 X_{\pi_1 \lambda_1}\cap g_2 X_{\pi_2 \lambda_2} \cap \dotsb \cap g_s
   X_{\pi_s \lambda_s} 
 \end{equation}
 is non-empty for any $g_1,\dotsc,g_s\in G$.
 Since the image in $G/P$ of this 
 intersection under the projection map 
 $\pr$~\eqref{E:FlagFibration} is a subset of 
 \begin{equation}\label{E:genint_G/P}
  g_1 \lambda_1^{-1}X_{\pi_1}\cap g_2 \lambda_2^{-1}X_{\pi_2} \cap \dotsb 
    \cap g_s \lambda_s^{-1}X_{\pi_s}\,,
 \end{equation}
 this latter intersection is non-empty for any $g_1,g_2,\dotsc,g_s\in G$,
 which proves (i).

 For (ii), let $p_1,\dotsc,p_s\in P$ be general.
 Then $p_1\lambda_1^{-1},\dotsc,p_s\lambda_s^{-1}$ are general, and the
 hypotheses imply that the intersection 
\[
  p_1 \lambda_1^{-1}X_{\pi_1}\cap p_2 \lambda_2^{-1}X_{\pi_2} \cap \dotsb 
    \cap p_s \lambda_s^{-1}X_{\pi_s}
\]  
 is transverse at the point $eP$.
 Similarly, the hypotheses imply that the intersection in $L/Q=P/R$
\[
  p_1 X_{\lambda_1}\cap p_2 X_{\lambda_2} \cap \dotsb \cap p_s X_{\lambda_s}
\]
 is non-empty and transverse at the generic point of each component.
 Thus Proposition~\ref{P:Transversality}(ii) implies that the intersection
\[
  p_1 X_{\pi_1\lambda_1}\cap p_2 X_{\pi_2\lambda_2} \cap \dotsb
   \cap p_s X_{\pi_s\lambda_s}
\]
 is transverse at a general point lying in the fiber $P/R$ above $eP$.
\end{proof}

Using Proposition~\ref{P:FeasibleSP}, we prove Theorem~\ref{T:necessary}.

\begin{proof}[Proof of Theorem~$\ref{T:necessary}$]
 By  Proposition~\ref{P:FeasibleSP}(ii), the Schubert position
 $\pi_1\lambda_1,\dotsc,\pi_s\lambda_s$ is feasible for $G/R$.
 Let $r_1,\dotsc,r_s$ be general elements of $R$.
 Then by Proposition~\ref{P:gen-feasible}
\[
   r_1 T_{\pi_1\lambda_1}\cap  r_2T_{\pi_2\lambda_2}
    \cap \dotsb \cap  r_s T_{\pi_2\lambda_s}
\]
is transverse.

 Since $\varphi_\frs$ is a surjection,
 Proposition~\ref{P:Transversality}(i)
 implies that the intersection
 \begin{equation}\label{E:general-w-transverse}
   \varphi_{\frs}(r_1 T_{\pi_1\lambda_1})\cap  
    \varphi_\frs(r_2T_{\pi_2\lambda_2})
    \cap \dotsb \cap  \varphi_\frs(r_s T_{\pi_s\lambda_s})
\end{equation}
 is transverse in $\frs^*$.
 This implies the codimension inequality 
\[
    \sum_{i=1}^s \codim\varphi_{\frs}(r_iT_{\pi_i\lambda_i})\ 
    \leq\ \dim \frs\,.
\]
 Since the map $\varphi_\frs$ is $R$-equivariant, these codimensions do
 not depend upon the choices of the $r_i$, which proves the theorem.
\end{proof}

\section{Proof of Theorem~\ref{Theorem}}\label{S:three}

\begin{figure}[htb]
 \begin{center}
   \begin{picture}(113,60)(-42,0)
    \put(0,0){\includegraphics{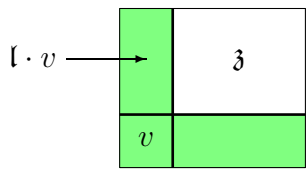}}
    \put(7,8){$v$}\put(43,39){$\frz$}
    \put(-42,38){$\frl\cdot v$}\put(-20,41){\vector(1,0){30}}
   \end{picture}
 \end{center}
 \caption{For the Grassmannian $\Gr(k,n)$, we may assume that $v$ is a rank $r$ 
 matrix concentrated in the lower left of 
 $T_{eP}\Gr(k,n)=\mbox{Mat}_{k\times(n{-}k)}$.
 Then $\frl\cdot v$ (shaded) and $\frz$ are as shown.}
 \label{F:ldotv}
\end{figure}
\begin{figure}[htb]
 \begin{center}
  \begin{picture}(165,130)(0,-20)
   \put(0,-20){\includegraphics{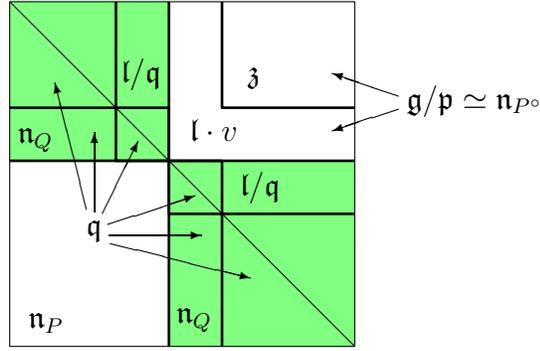}}
   \put(3, 57){$\frn_Q$}  \put(63,-11){$\frn_Q$}
   \put(7,-13){$\frn_P$}  \put(29,22){$\frq$}
    \put(29,32){\vector(-1,4){12}}
    \put(32,31){\vector(0,1){30}}
    \put(35,30.5){\vector(1,2){13.5}}
    \put(37,26){\vector(3,1){33}}
    \put(37.5,22){\vector(1,0){35}}
    \put(38,19){\vector(4,-1){53}}
   \put(42,80){$\frl/\frq$}  \put(87,37){$\frl/\frq$}
   \put(68,58){$\frl\cdot v$}
   \put(90,80){$\frz$}
   \put(150,70){$\frg/\frp\simeq\frn_{P^\circ}$}
   \put(147,75){\vector(-3,1){25}} \put(147,70){\vector(-3,-1){25}}
  \end{picture}
 \end{center}
 \caption{
 For the Grassmannian $\Gr(k,n)$, if $v$ has rank $r$, then 
 $L/Q\simeq\Gr(k{-}r,n)\times\Gr(n{-}k{-}r,n{-}k)$.
 We illustrate the weight decomposition of spaces
 $\frg/\frp$, $\frl \cdot v$, $\frz$, $\frq$, $\frl/\frq$, $\frn_P$, and 
 $\frn_Q$, where
 the off-diagonal entries in an $n \times n$ matrix represent
 the roots in $A_{n-1}$.  The roots of $\frl$ are shaded.}
\label{F:decompose}
\end{figure}
\begin{figure}[htb]
\begin{center}
   \begin{picture}(175,130)
    \put(0,0){\includegraphics{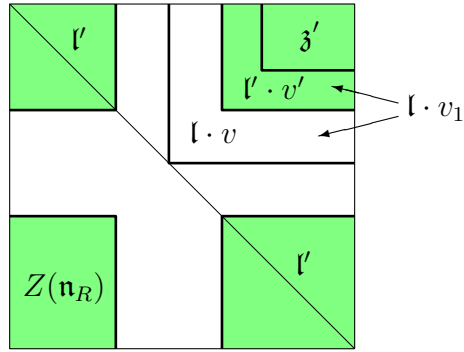}}
    \put(88,94.5){$\frl'\cdot v'$} 
    \put(110,115){$\frz'$} 
    \put(68,78){$\frl\cdot v$}
    \put(150,88){$\frl\cdot v_1$} 
            \put(147,93){\vector(-4,1){25}}
            \put(147,88){\vector(-4,-1){30}}
    \put(23,113){$\frl'$} 
    \put(108,28){$\frl'$} 
    \put(5,20){$Z(\frn_R)$}
   \end{picture}
\end{center}
\caption{
For the Grassmannian $\Gr(k,n)$, 
the semisimple part of $G'$ is isomorphic to 
$SL_{(k-r)+(n-k-r)}=SL_{n-2r}$,
whose roots are shaded.  
We also illustrate
the 
weights of $\frl'$, $\frl \cdot v$, $\frl' \cdot v'$, $\frz'$,
$\frl \cdot v_1$, and $Z(\frn_R)$.}
\label{F:Gprime}
\end{figure}

This proof is independent of Lie type and uses
some technical results
involving roots of the different groups ($G,P,R,L,Q,\dotsc$) and their 
Lie algebras, which we have collected together in the Appendix. 
For the classical groups, these results can
also be verified directly.  
For example, Lemma~\ref{L:Qcominuscule} shows that $L/Q$ is cominuscule
if $Q\in M(P)$;
this is also seen more concretely in Section~\ref{S:four} on a case-by-case
basis.
Figures~\ref{F:ldotv},~\ref{F:decompose}, and~\ref{F:Gprime}
illustrate the various groups and spaces that arise through
examples in type $A$.
In this case, $G/P=\Gr(k,n)$, the Grassmannian of $k$-planes in $\CC^n$,
the semisimple part of $L$ is $SL_k\times SL_{n-k}$, and the
tangent space at $eP$ is identified with $k\times(n{-}k)$ matrices.

We will prove Theorem~\ref{Theorem} in three stages, which we
formulate below.  

\begin{thm}\label{T:instages}
 Suppose that the semisimple part of $G$ is simple.   
 Let $\pi_1,\dotsc,\pi_s$ be a Schubert position for a cominuscule flag
 variety $G/P$.
 Then  $\pi_1,\dotsc,\pi_s$ is feasible if and
     only if any of the following equivalent conditions hold.
\begin{enumerate}
 \item[(i)] For every $Q \in M(P)\cup \{L\}$ and 
     every feasible Schubert position $\lambda_1,\dotsc,\lambda_s$ for $L/Q$,
     the intersection
 \begin{equation}\label{E:w-transverse}
   \varphi_R(r_1 T_{\pi_1\lambda_1})\cap  
    \varphi_R(r_2T_{\pi_2\lambda_2})
    \cap \dotsb \cap  \varphi_R(r_s T_{\pi_s\lambda_s})\,.
\end{equation}
  is transverse for general elements $r_1, \dotsc, r_s \in R$.

 \item[(ii)] For every  $Q \in M(P)\cup \{L\}$ and 
     every feasible Schubert position $\lambda_1,\dotsc,\lambda_s$ for $L/Q$,
     the inequality~\eqref{E:w-ineq} holds.

 \item[(iii)] For every $Q \in M(P) \cup \{L\}$ and 
     every feasible top-degree Schubert position
     $\lambda_1,\dotsc,\lambda_s$ for $L/Q$, 
     the inequality~\eqref{E:w-ineq} holds.
\end{enumerate}
\end{thm}

The intersection~\eqref{E:w-transverse} is the specialization 
of~\eqref{E:general-w-transverse} to the case where $\frs = Z(\frn_R)$,
and so the transversality of this intersection implies the 
inequality~\eqref{E:w-ineq}.
Thus the purely combinatorial statement of 
(ii) above
is {\em a priori} strictly stronger than (i), while (iii) is strictly 
stronger than (ii).  Theorem~\ref{T:instages}(iii) is precisely
Theorem~\ref{Theorem}.

Suppose that $\pi_1,\dotsc,\pi_s$ is a Schubert position for $G/P$ and 
$l_1,\dotsc,l_s$ are general elements of $L$.
By Proposition~\ref{Prop:feasible},  $\pi_1,\dotsc,\pi_s$ is feasible
if and only if the intersection
 \begin{equation}\label{E:TS-Int}
  T\ :=\ 
   l_1T_{\pi_1}\ \cap\  l_2T_{\pi_2}\ \cap\ \dotsb\ \cap\  l_sT_{\pi_s}
 \end{equation}
is transverse.

Since Theorem~\ref{T:necessary} establishes one direction of
Theorem~\ref{T:instages}, we assume that 
the Schubert position $\pi_1,\dotsc,\pi_s$ is
infeasible, and hence that the intersection~\eqref{E:TS-Int} is
non-transverse
when $l_1,\dotsc,l_s$ are general elements of $L$.
We first show that there is some
$Q \in M(P)$ and a feasible Schubert position
$\lambda_1,\dotsc,\lambda_s$ for $L/Q$ such that 
a general intersection~\eqref{E:w-transverse} is non-transverse.
This will prove Theorem~\ref{T:instages}(i).
Then, we use an inductive argument to show 
this implies that one of the inequalities~\eqref{E:w-ineq} is violated.

\subsection{A lemma on tangent spaces}

Since $L$ has only finitely many orbits on the tangent space
$\frg/\frp$, there is a unique largest orbit $O$
meeting the intersection $T$.
 This orbit does not depend on the generically chosen $l_1, \dotsc, l_s$.
Set  $V_i := (T_{\pi_i} \cap O)_\text{red}$ 
to be the variety underlying
the scheme-theoretic intersection of $T_{\pi_i}$ with this orbit.

 For any $v \in \frg/\frp$, we consider its $L$-orbit, $L\cdot v$.
As group schemes over $\CC$ are reduced, 
the tangent space to $L\cdot v$ at $v$ is $\frl\cdot v$.
Let $\frz$ be the quotient of $\frg/\frp$
by its subspace $\frl\cdot v$, and let 
$\psi\colon\frg/\frp\twoheadrightarrow\frz$ be the quotient map.

The main idea in our proof is the following result
concerning the images of the subspaces $l_iT_{\pi_i}$ in $\frz$.

\begin{lemma}\label{L:Inductive_engine}
Assume either that $v$ is a general point of\/ $T \cap O$,
or that $v$ is a smooth
point of each of the varieties $l_iV_i$.
 The intersection~\eqref{E:TS-Int} is transverse if and only if the
 intersection 
 \begin{equation}\label{E:image-int}
    \psi(l_1T_{\pi_1})\ \cap\  \psi(l_2T_{\pi_2})\ \cap\ \dotsb\ \cap\
    \psi(l_sT_{\pi_s})
 \end{equation}
 is transverse in the quotient space $\frz$.
\end{lemma}

Lemma~\ref{L:Inductive_engine} is invoked twice; once when $v$ is 
taken to be  a general point of $T \cap O$,
and a second time
when the varieties $l_iV_i$ are smooth at $v$ 
(but $v$ is chosen in advance, so {\em a priori} we do not know that it
is sufficiently general).
A consequence of our analysis is that smoothness of the 
$l_iV_i$ at $v$ is the condition for $v$ to be general.

We note that the intersection~\eqref{E:TS-Int} is transverse
if and only if for any $k \in L$, the intersection
\[
  kT\ =\ 
  (kl_1)T_{\pi_1}\ \cap\  (kl_2)T_{\pi_2}\ \cap\ \dotsb\ \cap\  (kl_s)T_{\pi_s}
\]
is transverse.  
When necessary we will therefore allow ourselves 
to translate $T$, and hence $v$, by an element of $L$.

\begin{remark}
Two special cases are worthy of immediate notice.  

Suppose that $v$ lies in
the dense orbit of $L$.  Then $\frz$ is zero-dimensional, and so
Lemma~\ref{L:Inductive_engine} implies that the 
intersection~\eqref{E:TS-Int} is necessarily transverse.

 On the other hand, suppose 
that $v=0$.  Then Lemma~\ref{L:Inductive_engine} 
provides no information.  However, since $v$ is assumed to lie in the 
largest orbit meeting $T$, we deduce that the subspaces
$l_iT_{\pi_i}$ meet only at the origin, and so 
$\sum \codim T_{\pi_i} \geq \dim \frg/\frp$.
Thus the intersection~\eqref{E:TS-Int} is transverse only when this is an
equality.  
\end{remark}

Since we assumed that the intersection~\eqref{E:TS-Int} is non-transverse,
we deduce that $v$ cannot lie in the dense orbit.
Moreover, if $v=0$, then the 
basic codimension inequality~\eqref{E:degreeineq}
arising from the degenerate case $Q=L$ is violated.  This second
observation will form the base case of the induction in our
proof of Theorem~\ref{T:instages}(ii).
Thus once we have proved Lemma~\ref{L:Inductive_engine} we will assume
that $v \neq 0$, and that $v$ does not lie in the dense orbit of $L$ on
$\frg/\frp$, as we have already dealt with these cases.

\subsection{Proof of Lemma~\ref{L:Inductive_engine}}
Under either hypothesis, we have $v \in O$, hence $O = L\cdot v$.
For each $i=1,\dotsc,s$, we consider the scheme-theoretic
intersection 
$l_iT_{\pi_i} \cap O$,
 whose underlying variety is 
$l_iV_i$.
Let $S_i$
denote the Zariski tangent space at $v$ to this scheme.
\[
   S_i\ :=\ T_v\bigl(l_iT_{\pi_i}\cap(L\cdot v)\bigr) \ =\ l_iT_{\pi_i}\cap (\frl\cdot v)
\]
Then $S_i\supset T_v(l_iV_i)$.

\begin{lemma}
  Under the hypotheses of Lemma~$\ref{L:Inductive_engine}$,
  the varieties $l_iV_i$ intersect transversally at $v$ in 
  $O$.
 Hence, the linear spaces $T_v(l_iV_i)$ are transverse in $\frl \cdot v$.
\end{lemma}

\begin{proof}
 Since $T \cap O$ is non-empty for generally chosen $l_1, \dotsc, l_s$,
 the intersection of general 
 $L$-translates of the varieties $V_i$ can never be empty.
 Since 
 $O$ is a homogeneous space of a reductive
 group, Kleiman's Transversality Theorem~\cite[Theorem  2(ii)]{Kl74} implies
 that the intersection of general $L$-translates 
 of the $V_i$ is transverse.
 The point $v$ lies in the intersection of the varieties
 $l_iV_i$.
 Since the elements $l_i\in L$ were chosen to be general,
 we conclude that the varieties $l_iV_i$ meet transversally at $v$,
 which by (either of)
 the hypotheses of Lemma~\ref{L:Inductive_engine} is a general 
 point of their intersection.
\end{proof}

\begin{cor}
   The linear subspaces $S_i$ are transverse in $\frl\cdot v$.
\end{cor}

Lemma~\ref{L:Inductive_engine} now follows from
Proposition~\ref{P:Transversality}(ii):
we have the exact sequence
\[
   0 \ \longrightarrow\ \frl \cdot v \ \longrightarrow\ 
   \frg/\frp\ \overset{\psi}{\longrightarrow}\ \frz \ \longrightarrow\ 0\,.
\]
with subspaces $l_iT_{\pi_i} \subset \frg/\frp$, and 
$S_i = l_iT_{\pi_i} \cap (\frl \cdot v)$ are transverse
 in $\frl \cdot v$.  
 \hfill\qed

\subsection{Proof of Theorem~\ref{T:instages}(i)}
We now show that Lemma~\ref{L:Inductive_engine} implies 
Theorem~\ref{T:instages}(i) by identifying the
intersection~\eqref{E:image-int} in $\frz$ with a general
intersection of the form~\eqref{E:w-transverse} in $Z(\frn_R)^*$,
for a parabolic subgroup $R$ of $P$ corresponding to 
some $Q \in M(P)$.

To this end, let $v$ be a general point of $T \cap O$, and
let $Q\subset L$ be the stabilizer of $\frl\cdot v$.
By Lemma~\ref{L:Qcominuscule}, $Q$ is a parabolic subgroup of $L$
and $L/Q$ is a cominuscule flag variety.
Translating $v$ by an element of $L$, we may furthermore assume that 
$Q$ is a standard parabolic, i.e. that $Q \supset B_L$.

Define $\lambda_i$ to be the Schubert position of $l_i^{-1}Ql_i$ with 
respect to $B_L$.
Then there exists a $b_i\in B_L$ such that 
$b_i^{-1}l_i^{-1}Ql_ib_i \supset \lambda_i B_L\lambda_i^{-1}$.
Set  $q_i := l_ib_i\lambda_i \in Q$.
Note that $\lambda_1,\dotsc,\lambda_s$ is automatically feasible, since
the $l_i$ are generic and $eQ$ lies in the intersection of the
translated Schubert cells $l_iB_L\lambda_iQ$.

%

By Corollary~\ref{C:frz} we have an $R$-equivariant isomorphism
$Z(\frn_R)^* \simeq \frz$.  

\begin{lemma}\label{L:Identify}
 We have 
 $\varphi_R(q_iT_{\pi_i\lambda_i})\simeq \psi(l_i T_{\pi_i})$.
\end{lemma}

\begin{proof}
Note that $B_L\subset \pi_i^{-1}B\pi_i$.  
Since $B_L\subset P$, it stabilizes both $P$ and 
$X_{\pi_i}$, and thus it stabilizes $T_{\pi_i}$.
We have the exact sequence~\eqref{Eq:Tan_SES} from 
Section~\ref{S:schubertvarieties},
\[
   T_{\lambda_i}\ \hookrightarrow\ T_{\pi_i\lambda_i}\ \twoheadrightarrow\
   \lambda_i^{-1}T_{\pi_i}\,. 
\]
Since $Q$ stabilizes the tangent spaces $\frl/\frq, \frg/\frr$, and $\frg/\frp$,
we may act on this sequence by $q_i:=l_ib_i\lambda_i$ to obtain
\[
   q_iT_{\lambda_i}\ \hookrightarrow\ q_iT_{\pi_i\lambda_i}\ 
    \twoheadrightarrow\ l_iT_{\pi_i}\,,
\]
as $b_i\in B_L$ stabilizes $T_\pi$.
This is a subdiagram of
 \[
  \begin{matrix}
    \frl/\frq & \into & \frg/\frr & \onto & \frg/\frp\\
        & & \makebox[0.01in][r]{\hspace{-15pt}\raisebox{2pt}{$\varphi_R$}}\Big\downarrow && 
            \Big\downarrow\makebox[0.01in][l]{\raisebox{2pt}{$\psi$}}\\
        & & Z(\frn)^*& \xrightarrow{\ \sim\ } & \frz
  \end{matrix}
 \]
We conclude that $\varphi_R(q_iT_{\pi_i\lambda_i})\simeq \psi(l_i T_{\pi_i})$, 
under the identification
of $Z(\frn_R)^*$ with $\frz$. 
\end{proof}

Since the intersection~\eqref{E:TS-Int} is assumed to be non-transverse,
Lemma~\ref{L:Inductive_engine} implies that 
the intersection~\eqref{E:image-int} is non-transverse.
Lemma~\ref{L:Identify} shows that this is equivalent to
\[
  \varphi_R(q_1T_{\pi_1\lambda_1})\ \cap\ 
  \varphi_R(q_2T_{\pi_2\lambda_2})\ \cap\ \dotsb\ \cap
  \varphi_R(q_sT_{\pi_s\lambda_s})
\] 
being non-transverse.

This is an intersection of the form~\eqref{E:w-transverse}, however,
since the $q_i$ are constructed from $v$ and $l_i$, they will not be
general elements of $R$ (they are not even general elements of $Q$).  
It remains to show
that a {\em general} intersection~\eqref{E:w-transverse} is 
non-transverse.

Consider what happens when we translate each $l_i$ by a general
element $k_i \in \Stab_L(\CC v) \subset Q$.  
The point $v$ will still be a point
of the new intersection
\[
   T'\ :=\ (k_1l_1)T_{\pi_1}\ \cap\  (k_2l_2)T_{\pi_2}\ \cap\ \dotsb\ \cap\
   (k_sl_s)T_{\pi_s} \,,
\]
thus we obtain the same subgroup $Q$.
Moreover, since $v$ is a smooth point of $l_iV_i$, 
and the $k_i$ are general,
it will be a smooth point
of $(k_il_i)V_i$.
If $q_i'$ denotes the new $q_i$ we obtain for the intersection $T'$, 
we find that $q_i' = k_i q_i$.
Thus by Lemmas~\ref{L:Inductive_engine} and~\ref{L:Identify} we see 
that the intersection
\[
   \varphi_R(k_1q_1T_{\pi_1\lambda_1})\ \cap\
   \varphi_R(k_2q_2T_{\pi_2\lambda_2})\ \cap\ \dotsb\ \cap\  
   \varphi_R(k_sq_sT_{\pi_s\lambda_s})
\]
is non-transverse for general $k_i \in \Stab_L(\CC v)$.
By Lemma~\ref{L:aseffective}, this implies that a 
general intersection~\eqref{E:w-transverse} is non-transverse.
This proves Theorem~\ref{T:instages}(i).
\hfill\qed


\subsection{Proof of Theorem~\ref{T:instages}(ii)}

Recall that
$M(P)$ is exactly the set of those standard parabolic subgroups
of the form $\Stab_L(\frl \cdot v)$ 
for some $v \in \frg/\frp$.

We show that if $\pi_1,\dotsc,\pi_s$ is an infeasible 
Schubert position for $G/P$, then there is a parabolic subgroup $Q\in M(P)$ of
$L$ and a feasible Schubert position $\lambda_1,\dotsc,\lambda_s$ for $L/Q$ such that
 \begin{equation}\label{Eq:ineq}
   \sum_{i=1}^s \codim \varphi_R(T_{\pi_i \lambda_i})\ >\ 
    \dim Z(\frn_R)\,,
 \end{equation}
where $R$ is the parabolic subgroup of $P$ containing $Q$.

Suppose that Theorem~\ref{T:instages}(ii) holds for any proper
subgroup of $G$ whose semisimple part is simple,
and let $\pi_1,\dotsc,\pi_s$ be an infeasible Schubert position for
$G/P$. 
By  Theorem~\ref{T:instages}(i), there is a parabolic subgroup $Q \in M(P)$ and
a feasible Schubert position $\lambda_1,\dotsc,\lambda_s$ for $L/Q$ such that
for general  $r_1,\dotsc,r_s\in R$,
the intersection
 \begin{equation}\label{E:varphi}
   \varphi_R(r_1 T_{\pi_1\lambda_1})\cap  \varphi_R(r_2T_{\pi_2\lambda_2})
    \cap \dotsb \cap  \varphi_R(r_s T_{\pi_s\lambda_s})
 \end{equation}
is not transverse.
If this intersection has dimension $0$, then we deduce the codimension
inequality~\eqref{Eq:ineq}
and so we are done.

Now we assume that the dimension of the intersection~\eqref{E:varphi} is not
zero, and we use our inductive hypothesis to find a different parabolic 
subgroup $Q_1\in M(P)$ and a feasible Schubert position 
$\mu_1, \dotsc, \mu_s$ for $L/Q_1$ so that the corresponding inequality holds.  

We begin by constructing a new cominuscule flag variety $G'/P'$
whose tangent space at $eP'$ is identified with $\frz$.  This will
allow us to identify the intersection~\eqref{E:varphi} as an
intersection of tangent spaces of Schubert varieties.
Define the reductive (proper) subgroup $G'$ of $G$ to be
\[
   G'\ :=\ Z_G ( Z_H ( Z(N_R) ))\,.
\]
 $G'$ is the smallest reductive subgroup of $G$ containing both
 $H$ and $Z(N_R)$.
 Set $P':= G'\cap R$.  
Let $L'$ denote the Levi subgroup of $P'$,
and let $W'$ denote the Weyl group of $G'$.

By Lemma~\ref{L:sspart} the semisimple part of $G'$ is simple and
by Lemma~\ref{L:nilr} $G'/P'$ is cominuscule.
Thus the inductive hypothesis applies to $G'/P'$.

The pattern map  $w\mapsto \overline{w}$ of Billey and
Braden~\cite{BB02} sends $W\to W'$.
The element $\overline{w}\in W'$ is defined by its inversion set, which is 
$\Phi(\frg')\cap \Inv(w)$.

\begin{lemma}\label{L:tangent_equal}
   For all $w\in W^R$, $\varphi_R(T_w)=T_{\overline{w}}$. 
\end{lemma}

\begin{proof}
 Since $w\in W^R$, $\Inv(w)=\Phi(T_w)$.
 The weights of the tangent space $T_{\overline{w}}$ are the inversions of 
 $\overline{w}$ which lie in $\Phi(\frg'/\frp')$.
 By Lemma~\ref{L:nilr}, $\Phi(\frg'/\frp')=\Phi(\frz)$.
 Since the weights of $\varphi_R(T_w)$ are $\Inv(w)\cap\Phi(\frz)$, we 
  are done.
\end{proof}

By Lemmas~\ref{L:tangent_equal} and~\ref{L:aseffective} there exist 
$l'_1,\dotsc,l'_s\in L'$ such that 
the intersection~\eqref{E:varphi} is equal to
\[
   l'_1 T_{\overline{\pi_1\lambda_1}}\cap  
   l'_2 T_{\overline{\pi_2\lambda_2}} \cap \dotsb \cap  
   l'_s T_{\overline{\pi_s\lambda_s}}\,.
\]
Furthermore, as the elements $r_i\in R$ are general, so are the elements
$l'_i\in L'$.
Since this intersection is not transverse, we conclude that
if we set $\pi'_i:=\overline{\pi_i\lambda_i}$, then 
$\pi'_1,\dotsc,\pi'_s$ is an infeasible Schubert position for $G'/P'$.

By our inductive hypothesis, 
there is a parabolic subgroup $Q'\in M(P')$ and
feasible Schubert positions $\lambda'_1,\dotsc,\lambda'_s$ 
such that
 \begin{equation}\label{E:inductive_ineq}
   \sum_{i=1}^s \codim \varphi'_{R'}(T_{\pi'_i\lambda'_i})
    \ >\ \dim Z(\frn_{R'})\,.
 \end{equation}
(Here, $R'\subset P'$ is the largest parabolic subgroup such that 
$R'\cap L'=Q'$.)
Then $Q'$ is a standard parabolic which stabilizes $\frl'\cdot v'$ 
for some $v' \in \frg'/\frp' (\simeq \frz)$.

Let $Q_1$ be the stabilizer in $L$ of $\frl\cdot v_1$, where 
$v_1=v+v'$ (we consider $v'$ to be an element of $\frg/\frp$ by the
$L'$-equivariant injection $\frg'/\frp' \into \frg/\frp$).
It follows from Lemma~\ref{L:Q1standard} that $Q_1$ is a standard
parabolic, and so $Q_1\in M(P)$.  Let $R_1$ be the corresponding
parabolic subgroup of $P$.  By Lemma~\ref{L:indfrn},
$Z(\frn_{R_1})=Z(\frn_{R'})$, and 
$\frz' = (\frg'/\frp')/(\frl' \cdot v')$ is the dual to this space.

Let $\mu_i$ be the minimal coset representative of $\lambda_i\lambda'_i$
in $W_L/W_{Q_1}$.  Since 
$\lambda_1, \dotsc, \lambda_s$ is 
feasible for $L/Q = P/R$, and 
$\lambda'_1, \dotsc, \lambda'_s$ is feasible for
$L'/Q' = R/(R\cap R_1)$, 
$\lambda_1\lambda'_1, \dotsc, \lambda_s\lambda'_s$ is feasible
for $P/(R \cap R_1)$, by Proposition~\ref{P:FeasibleSP}(ii).
Hence by Proposition~\ref{P:FeasibleSP}(i),
$\mu_1, \dotsc, \mu_s$ is feasible for $P/R_1 = L/Q_1$.

We now complete the proof by showing that 
$\dim \varphi_{R_1}(T_{\pi_i \mu_i})
 =\dim \varphi'_{R'}(T_{\pi'_i \lambda'_i})$.
These $H$-invariant subspaces have weights 
$\Inv(\pi_i \mu_i)\cap \Phi(\frz')$ and 
$\Inv(\overline{\pi_i \lambda_i} \lambda'_i)\cap \Phi(\frz')$,
respectively.
Let $\nu_i = \mu_i^{-1}\lambda_i\lambda_i' \in W_{Q_1}$.  Then
by~\eqref{E:inversionequation},
\[
   \Inv(\pi_i \mu_i) \cap \Phi(\frz') \ 
   =\ \bigl(\nu_i \Inv(\pi_i \lambda_i \lambda'_i)\bigr) \cap \Phi(\frz') \ 
   =\ \nu_i \bigl(\Inv(\pi_i \lambda_i \lambda'_i) \cap \Phi(\frz')\bigr)\,,
\]
as $W_{Q_1}$ preserves $\Phi(\frz')$.  Thus it suffices
to show that 
\[
  \Inv(\pi_i \lambda_i \lambda'_i) \cap \Phi(\frz')\  =\ 
   \Inv(\overline{\pi_i \lambda_i}\, \lambda'_i)\cap \Phi(\frz')\,.
\]
Note that we have 
$\overline{\pi\lambda \lambda'}=\overline{\pi\lambda} \lambda'$, as the
pattern map is $W'$-equivariant.  Then indeed
 \begin{align*}
  \Inv(\pi\lambda \lambda')\cap \Phi(\frz')\ 
    &=\ 
     (\pi\lambda \lambda')^{-1}\Phi^-  \cap \Phi^+(\frg') \cap \Phi(\frz')\\
    &= \ 
     \Inv(\overline{\pi\lambda \lambda'}) \cap \Phi(\frz')\\
    &=\ 
      \Inv(\overline{\pi\lambda} \lambda')\cap \Phi(\frz')\,.
\end{align*}
Thus we have exhibited a parabolic subgroup $Q_1 \in M(P)$ and a feasible
Schubert position $\mu_1, \dotsc, \mu_s$ for $L/Q_1$, such
that by rewriting~\eqref{E:inductive_ineq} we have
\[
  \sum_{i=1}^s \codim \varphi_{R_1}(T_{\pi_i\mu_i})
  \ >\  \dim Z(\frn_{R_1})\,,
\]
as required.  
\hfill\qed

\subsection{Proof of Theorem~\ref{T:instages}(iii)}

We need the following non-obvious fact which is proven in the
Ph.D. Thesis~\cite{KP_Thesis}.

\begin{prop}\label{P:KP_thesis}
  Suppose that $\pi'<\pi$ in the Bruhat order.
  Then there is an injection $\iota \colon \Inv(\pi')\into\Inv(\pi)$
  such that if   $\alpha\in\Inv(\pi')$, then $\iota(\alpha)$ is a higher root
  than $\alpha$.
\end{prop}
 
\begin{proof}[Sketch of Proof]
 It is enough to show this when $\pi'$ covers $\pi$ in the Bruhat order.
 In this case, $\pi'$ and $\pi$ differ by reflection in a root
 $\beta$, and one can verify the proposition by comparing inversions
 within strings of roots along lines parallel to $\beta$.  
\end{proof}

 Let $\pi_1,\dotsc,\pi_s$ be an infeasible Schubert position for $G/P$.
 Then by Theorem~\ref{T:instages}(ii), there exists a parabolic subgroup $Q\in
 M(L)$ and a feasible Schubert position $\lambda_1,\dotsc,\lambda_s$ for $L/Q$
 such that the inequality~\eqref{Eq:ineq} holds.

 If this Schubert position for $L/Q$ does not have top-degree, then 
 by Chevalley's formula~\cite{Ch91}, there exists a feasible Schubert position
 $\mu_1,\dotsc,\mu_s$ for $L/Q$ such that $\mu_i\leq \lambda_i$, for 
 $i=1,\dotsc,s$.
 Since each $\pi_i$ is a minimal coset representative, we have
 $\pi_i \mu_i\leq \pi_i \lambda_i$.
 Recall that the dimension of $\varphi_R(T_{\pi_i\lambda_i})$ is the number of
 inversions of $\pi_i\lambda_i$ which lie in the set of weights $\Phi(\frz)$.
 Since $N_R$ is $B$-stable, so is its center $Z(N_R)$,
 and hence the roots in $\Phi(\frz)=-\Phi(Z(\frn_R))$ are an upper order 
 ideal in $\Phi(\frg)$.
 Then Proposition~\ref{P:KP_thesis} implies that 
 $\dim \varphi_R(T_{\pi_i\mu_i}) \leq \dim \varphi_R(T_{\pi_i\lambda_i})$,
 and thus~\eqref{Eq:ineq} holds for $\mu_1,\dotsc,\mu_s$ in place of 
 $\lambda_1,\dotsc,\lambda_s$.
 \hfill\qed

%
\section{Explicating the Horn recursion}\label{S:four}
By Theorem~\ref{Theorem}, the feasibility of a Schubert position 
$\pi_1,\dotsc,\pi_s$ for cominuscule
$G/P$ is detected by the inequality~\eqref{E:w-ineq} for
every feasible top-degree Schubert position $\lambda_1,\dotsc,\lambda_s$ for $L/Q$ for
every $Q\in M(P)$. 
We noted in Remark~\ref{R:combinatorial} that  
these inequalities are combinatorial conditions.  
We now reformulate this. 
Write $\Invc(\pi)$ for the set of weights $\Phi(\frg/\frp)-\Inv(\pi)$
and call these the \Blue{{\em coinversions}} of $\pi$.
They are the weights of the normal bundle, $(\frg/\frp)/T_\pi$, to
$X_\pi$ at $eP$. 

\begin{lemma}\label{L:reformulate}
 Given a Schubert position $\pi_1,\dotsc,\pi_s$ for $G/P$ and a feasible
 Schubert position  $\lambda_1,\dotsc,\lambda_s$ for $L/Q$ with $Q\in M(P)$, the
 inequality~\eqref{E:w-ineq}  is equivalent to 
 \begin{equation}\label{Eq:New_Horn}
    \sum_{i=1}^s 
     \bigl| \Invc(\pi_i)\cap \lambda_i\Phi(\frz) \bigr|\ \leq\ \dim \frz\,,
 \end{equation}
 where $\frz = Z(\frn_R)^*$.
\end{lemma}

\begin{proof}
 As we observed in Remark~\ref{R:combinatorial}, the 
 inequality~\eqref{E:generalineq}
 (and hence~\eqref{E:w-ineq})
 can be computed combinatorially 
 as
 $\codim \varphi_\frs(T_{\pi_i\lambda_i}) = 
 |\Phi(\frs^*) - \Inv(\pi_i\lambda_i)|$.  
 Since 
 $\frs^* = \frz$, 
 by~\eqref{E:inversionequation} 
we have 
 \begin{align*}
  \codim \varphi_R(T_{\pi_i\lambda_i})\ 
   & =\  |\Phi(\frz) - \Inv(\pi_i\lambda_i)| \\
   & =\  |\Phi(\frz) \cap (\Phi(\frg/\frp) - \Inv(\pi_i\lambda_i))| \\
   & =\  |\Phi(\frz) - \lambda_i^{-1} \Inv(\pi_i)|  \\
   & =\  |\Phi(\frz) \cap \lambda_i^{-1} \Invc(\pi_i)|\ .
 \end{align*}
Translating by $\lambda_i$, this is equal to 
  $|\lambda_i \Phi(\frz) \cap \Invc(\pi_i)|$, which implies the lemma.
\end{proof}

We introduce the following notation.
Given a Schubert position $\pi$ for $G/P$ and a Schubert position $\lambda$ for
$L/Q$, 
set $|\pi|_\lambda := |\Invc(\pi)\cap \lambda\Phi(\frz)|$.
We also write 
$|\pi| := |\Invc(\pi)| = \codim T_\pi$.
Then the inequalities of Lemma~\ref{L:reformulate} become
\[
    \sum_{i=1}^s |\pi_i|_{\lambda_i}\ \leq\ \dim \frz\,.
\]
whereas the basic codimension inequality~\eqref{E:degreeineq} becomes
\[
    \sum_{i=1}^s |\pi_i|\ \leq\ \dim \frg/\frp\,.
\]

Since $G/P$ is cominuscule, the weights $\Phi(\frg/\frp)$ form a
lattice~\cite{P84}. 
For $\pi\in W^P$, the tangent space $T_\pi$ is $B_L$-invariant, so its weights
form a lower order ideal in this lattice.
Given a poset $Y$, let $J(Y)$ be the distributive lattice of lower order ideals
of $Y$~\cite{St86}.
Proctor~\cite{P84} showed that 

\begin{prop}\label{P:Order_Ideals}
 $W^P\simeq J(\Phi(\frg/\frp))$.
\end{prop}

\begin{remark}
 Proposition~\ref{P:Order_Ideals} allows us to interpret the
 inequalities~\eqref{E:w-ineq} in terms of convex geometry.
 Let $V$ be the vector space of functions $f\colon\Phi(\frg/\frp)\to\RR$.
 The set 
\[
   \calO_{\frg/\frp}\ :=\ 
  \{ f\in V\mid \alpha<\beta\in\Phi(\frg/\frp)\ 
     \Rightarrow\ 0\leq f(\alpha)\leq f(\beta)\leq 1\}
\]
 of order preserving maps from $\Phi(\frg/\frp)$  to $[0,1]$ 
 is the order polytope~\cite{St86b} of the poset $\Phi(\frg/\frp)$.
 Its integer points are the indicator functions of upper order ideals in 
 $\Phi(\frg/\frp)$, which by Proposition~\ref{P:Order_Ideals} are the
 indicator functions of the coinversion sets $\Invc(\pi)$ of Schubert
 positions $\pi$ for $G/P$.
 Write $u_\pi\in V$ for the integer point of $\calO_{\frg/\frp}$
 corresponding to the Schubert position $\pi$.

 Given a Schubert position $\lambda$ for $L/Q$ with $Q\in M(P)$, 
 define a linear map $\Sigma_\lambda\colon V \to \RR$
 by 
\[
  \Sigma_\lambda( f )\ :=\ \sum_{\gamma\in \lambda\Phi(\frz)} f(\gamma)\ .
\]
 Then $|\pi|_\lambda=\Sigma_\lambda(u_\pi)$.

 In particular, the inequality~\eqref{E:w-ineq} may be interpreted as a
 linear inequality on the polytope $(\calO_{\frg/\frp})^s$,
 and so the set of all feasible Schubert positions
 $\pi_1,\dotsc,\pi_s$ for $G/P$ is naturally identified with the
 integer points in the \Blue{{\em feasibility polytope}} which is the
 subpolytope of $(\calO_{\frg/\frp})^s$ defined by the set of
 inequalities from Theorem~\ref{Theorem}.
 We have not studied the structure of this feasibility polytope.

\end{remark}

We now investigate the inequalities of Theorem~\ref{Theorem} on a 
case-by-case basis.  
Recall that $M(P)$ is the set of standard
parabolic subgroups of $L$
of the form $Q = \Stab_L(T_v L\cdot v)$, for some $v \in \frg/\frp$.  Any two
suitable choices of $v$ in the same $L$-orbit give the same $Q$.  Thus 
for each type, it is enough analyze one such choice of $v$ from each 
$L$-orbit.  
The cases where $v=0$ or $v$ is in the dense orbit can be excluded,
since these yield $\Stab_L(T_v L\cdot v) = L$.
We can always take $v$ to be of the form
\[
  v\ =\ v_{\alpha_1}+\dotsb+v_{\alpha_r}\,,
\] 
where $v_\alpha \in \frg/\frp$ is a non-zero vector of weight $\alpha$,
and $\alpha_1, \dotsc, \alpha_r$ is a sequence of orthogonal long roots.  
The number $r$ determines the $L$ orbit of $v$~\cite{RRS92}.  We will
also make use of Lemma~\ref{L:frlcdotv}, which asserts that for such
a choice $v$, the weights of $\frz$ will be the weights of $\frg/\frp$
orthogonal to $\alpha_1, \dotsc, \alpha_r$.

\subsection{Type $A_{n-1}$, the classical Grassmannian, $\Gr(k,n)$.}\label{S:An}
Suppose that $P$ is obtained by omitting the $k$th node in the Dynkin diagram of
$A_{n-1}$.
Then $G/P$ is $\Gr(k,n)$, the Grassmannian of $k$-planes in $\CC^n$. 
The Levi subgroup $L$ of $P$ 
has semisimple part $SL_k\times SL_{n-k}$.
We identify $\frg/\frp$ with $\Hom(\CC^k,\CC^{n-k})$,
where $\CC^k\oplus\CC^{n-k}=\CC^n$.
Its weights are 
\[
   \Phi(\frg/\frp)\ =\ 
   \{ e_j-e_i \mid  1\leq i\leq k < j \leq n \}\,,
\]
where $e_1,\dotsc,e_n$ are the standard orthonormal basis vectors of
$\CC^n=\frh^*$.
We identify $\Phi(\frg/\frp)$ with the cells of a 
$k\times(n{-}k)$ rectangle where $e_j-e_i$ corresponds to the cell in
row $i$ (from the top) and column $j{-}k$ (from the left).
The lowest root in 
$\Phi(\frg/\frp)$
is in the lower left corner and the highest root is in
the upper right corner.

Minimal coset representatives $\pi\in W^P$ are permutations of $n$ with
a unique descent at position $k$.
The inversion set of a permutation $\pi$ is the set of roots
\[
   \{ e_j-e_i \mid i\leq k<j\mbox{\ such that\ }\pi(i)>\pi(j)\}\,.
\]
We display this for $n=11$, $k=5$, 
and $\pi=1367\,\ten\,24589\,\eleven$, 
shading the inversion set.
 \begin{equation}\label{Eq:1367}
  \raisebox{-25pt}{%
  \begin{picture}(177,60)
   \put(48,0){\includegraphics{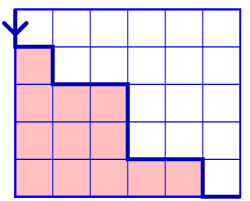}}
   \put(  0,14){$\Inv(\pi)$}   \put( 35,17){\vector(1,0){25}}
   \put(138.5,26){$\Invc(\pi)$} \put(135.5,29){\vector(-1,0){25}}
  \end{picture}}
 \end{equation}
The permutation may be read off from the inversion diagram as follows.
Consider the path which forms the border of $\Inv(\pi)$ from the upper
left corner to the lower right corner of the rectangle.
If we label the steps from 1 to $n$, then the labels of the vertical
steps are the first $k$ values of $\pi$ and the labels of the horizontal
steps are the last $n{-}k$ values of $\pi$.

 If we write $\alpha_{i}:=e_{k+i}-e_{k+1-i}$,
 which is the $i$th 
 root along the the anti-diagonal 
 in $\Phi(\frg/\frp)$ starting from the lower left,
 then the vector $v$ may be taken to have the form
 \begin{equation}\label{Eq:vee}
  v\ =\ 
       v_{\alpha_{1}}+ v_{\alpha_{2}}+ \dotsb  +v_{\alpha_{r}}\,,
 \end{equation}
and  $L\cdot v\subset \Hom(\CC^k,\CC^{n-k})$ consists of rank $r$
matrices.
Note that $1 \leq r < \min\{k,n{-}k\}$.
Then the set $\Phi(\frz)$ is the upper right
$(k{-}r)\times(n{-}k{-}r)$ rectangle in the $k\times(n{-}k)$ rectangle
representing $\Phi(\frg/\frp)$, and $\dim \frz=(k{-}r)(n{-}k{-}r)$. 
We show this for $n=11$, $k=5$, and $r=2$.
\[
  \begin{picture}(272,58)(-108,0)
   \put(0,0){\includegraphics{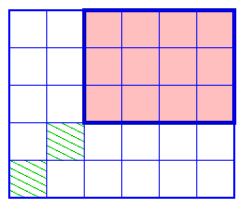}}
   \put(86.5,36){$\Phi(\frz)$} \put(83.5,39){\vector(-1,0){25}}
   \put(86.5, 2){$\Phi(\frs\frl_k\times\frs\frl_{n-k}\cdot v)$} 
   \put(83.5, 6){\vector(-1,0){25}}
   \put(-36,15){$\alpha_{2}$} \put(-22,17){\vector(1,0){41}}
   \put(-36, 4){$\alpha_{1}$} \put(-22, 6){\vector(1,0){30}}
  \end{picture}
\]

The subgroup $Q\in M(P)$ 
which is the stabilizer of 
$\frs\frl_k\times\frs\frl_{n-k}\cdot v$ is obtained by 
further omitting the nodes at $k{-}r$ and at $k{+}r$ in the Dynkin
diagram for $L^{ss}$.
\[
  \begin{picture}(276,43)(-35,-20)
   \put(-36.5,0){\includegraphics{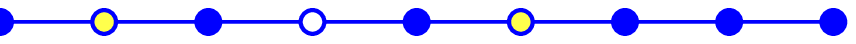}}
   \put( 86,14){$P$}   \put(85,-20){$Q$}
   \put( 97,-14){\vector( 3,1){46}}  
   \put( 81,-14){\vector(-3,1){46}}
  \end{picture}
\]
Thus $L/Q$ is isomorphic to $\Gr(k{-}r,k)\times\Gr(r,n{-}k)$.

An element $\lambda\in W^Q$ acts on $\Phi(\frg/\frp)$ by simultaneously
shuffling the $r$ rows 
that do not meet $\Phi(\frz)$ with those that do, and
the same for columns. 
This is equivalent to selecting $r$ rows and $r$ columns, the images under
$\lambda$ of the rows and columns which do not meet $\Phi(\frz)$.
If we draw 
$\Invc(\pi)$ in the rectangle and  cross out the selected rows and columns,
then $|\pi|_\lambda$ is the number of boxes which
remain. 
In the example~\eqref{Eq:1367} above with $\pi=1367\,\ten\,24589\,\eleven$ 
and $r=2$, if  $\lambda$ selects rows 2 and 4 
from the top and columns 2 and 6
from the right, we see that $|\pi|_\lambda=7$.
\[
  \begin{picture}(80,68)(-10,-8)
    \put(0,0){\includegraphics{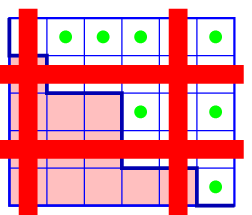}}
    \put(-10,37){$2$}  \put( 5,-10){$6$}
    \put(-10,16){$4$}  \put(49,-10){$2$}
  \end{picture}   
\]

\begin{remark}\label{Rem:rem}
For the purpose of our cominuscule recursion we describe how to obtain the
inversion diagram of an element $\lambda\in W^Q$, which is a subset of a
$(k{-}r)\times r$ 
rectangle for the rows and a $r\times(n-k-r)$ rectangle for the
columns. 
In the rectangle for the rows, draw a path from the 
upper left corner to the lower right corner 
whose $i$th step is {\em horizontal} 
if $\lambda$ selected row $i$
and vertical otherwise, while in the rectangle for the columns, 
draw 
a path from the lower right corner to the upper left corner 
whose $i$th step is {\em vertical} if $\lambda$ selected 
column $i$ and horizontal otherwise.
We show this for our example.
\[
  \begin{picture}(129,47)
   \put(0,0){\includegraphics{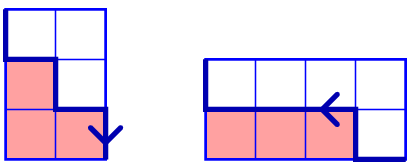}}
    \put(  5,34){2} \put( 20,20){4}
    \put( 63,20){6} \put(105, 6){2} 
  \end{picture}
\]
\end{remark}

Since $L/Q$ is a product of smaller cominuscule flag varieties, feasibility 
for the Schubert positions $\lambda$ in the recursion is determined separately
on each factor.  
Note that not all cominuscule $L/Q$ enter into this recursion.

\subsection{Type $D_{n+1}$, $G = SO_{2n+2}$, $G/P$ is the 
even-dimensional quadric, $Q^{2n}$.}\label{S:even-quadric}
  
Here, the parabolic subgroup $P$ is obtained by omitting the rightmost
node of the Dynkin diagram, as shown in Table~\ref{Ta:One}.
Its Levi subgroup $L$ has semisimple part $SO_{2n}$ and the flag variety $G/P$ is the
even-dimensional quadric $Q^{2n}$ in $\PP^{2n+1}$.
The lattice $\Phi(\frg/\frp)$ is the poset $\Lambda_{n-1}$, 
whose Hasse diagram we display,
where elements to the right are greater.
\[
  \begin{picture}(274,55)(-2,-25)
   \put( 0,0){\circle*{5}} \put(30,0){\circle*{5}}
   \put( 0,0){\line(1,0){30}} \put(30,0){\line(1,0){20}}
   \put(59,-1){$\ldots$}
   \put(80,0){\line(1,0){20}} \put(160,0){\line(1,0){20}}
   \put(100,0){\circle*{5}} \put(130, 10){\circle*{5}}
                            \put(130,-10){\circle*{5}} \put(160, 0){\circle*{5}}
   \put(100,  0){\line(3,1){30}} \put(100, 0){\line(3,-1){30}}
   \put(130,-10){\line(3,1){30}} \put(130,10){\line(3,-1){30}}
   \put(189,-1){$\ldots$}

   \put(230,0){\circle*{5}} \put(260,0){\circle*{5}}
   \put(230,0){\line(1,0){30}}\put(210,0){\line(1,0){20}}

   \put(-3,-15){$0$}   \put(27,-15){$1$}  \put(245,-15){$2n{-}2$}
   \put(119,-24){$n{-}1$}  \put(119, 17){$\overline{n{-}1}$}
   \put(157,-15){$n$} 
  \end{picture}
\]
Each root %
in $\Phi(\frg/\frp)$
is orthogonal to exactly one other, and their indices sum to $2n{-}2$.
Consequently an orthogonal sequence of long roots has length at most 2.
For our purposes, there is one interesting orbit of $L$ in $\frg/\frp$.
In fact, $\frg/\frp$  is the defining representation of $SO_{2n}$ and this orbit is the
set of (non-zero) isotropic vectors, the cone over the quadric $Q^{2n-2}$.
Thus $M(P)$ consists of a single parabolic subgroup $Q$, where  $L/Q$ is
the quadric $Q^{2n-2}$,  and $W^Q=\Lambda_{n-1}$.  Here $Q$ is
the stabilizer in $L$ of $\frl \cdot v_\alpha$, where $\alpha$ is
the simple root defining $P$ 
(labeled $0$ in $\Lambda_{n-1}$) 
and $\Phi(\frz)$ is the orthogonal complement to $\alpha$ which is the 
single root labelled $2n{-}2$ in $\Lambda_{n-1}$. 

By Proposition~\ref{P:Order_Ideals}, 
$W^P$ is the set of order
ideals of $\Lambda_{n-1}$, which is equal to
$\Lambda_n$, where the set of weights of $T_\pi$ is equal to the order 
ideal $\pi$, and $|\pi|$ is the cardinality of the {\em complement} of
this order ideal.  Thus $\lambda \in W^Q$ is an element of $\Lambda_{n-1}$, 
whereas $\pi \in W^P$ is an order ideal of $\Lambda_{n-1}$.  The action
of $W_P$ on both $\Phi(\frg/\frp)$ and $W_Q$ canonically identifies
these two occurrences of $\Lambda_{n-1}$; however, as the identification
of $W^P$ with $\Lambda_n$ is not canonical, there is a choice to be made.
We will adopt the convention that $n \in \Lambda_n$ corresponds to
the $n$-element order ideal in $\Lambda_{n-1}$ which contains $n{-}1$ 
and $\overline{n}$ corresponds to the $n$-element order ideal which contains 
$\overline{n{-}1}$.
For $\lambda\in W^Q=\Lambda_n$, we see that $\lambda\Phi(\frz)$ is the
root $\lambda^\perp$ orthogonal to $\lambda$, which is found by rotating
the Hasse diagram by $180^\circ$.
Thus
 \begin{equation}\label{E:even-quadric}
   |\pi|_\lambda\ =\ 
   \left\{\begin{array}{rcl}
     0&\ &\mbox{if } \lambda^\perp \in \pi\\
     1&&\mbox{otherwise}
\end{array}\right.
 \end{equation}
%
For example,
 \begin{equation}\label{E:hats}
   |n|_{\overline{n{-}1}}\ =\ |\overline n|_{n-1}\ =\ 1 
    \quad \text{and} \quad
   |n|_{n-1}\ =\ |\overline n|_{\overline{n{-}1}}\ =\ 0\,.
 \end{equation} 

Since $|M(P)|=1$ and $L/Q$ is $Q^{2n-2}$, the cominuscule recursion in 
this case can
proceed by induction on $n$.
The base case is $Q^2$, the quadric in $\PP^3$ which is isomorphic to
$\PP^1\times\PP^1$. 

Note that the condition 
\[
\sum_{i=1}^s |\pi_i|_{\lambda_i}\ \leq\ 1\,,
\] 
for $\lambda_1, \dotsc, \lambda_s$ feasible for $L/Q$, is implied by
the basic codimension inequality 
\[
  \sum_{i=1}^s |\pi_i|\ \leq\ 2n\,,
\]
unless $|\pi_1| + |\pi_2| = 2n$ and $|\pi_3| = \dotsb = 0$ 
(or some permutation thereof).  Indeed if 
$\lambda_1, \dotsc, \lambda_s$ is feasible for $L/Q$, and 
$|\pi_1|_{\lambda_1} = |\pi_2|_{\lambda_2} = 1$, then 
$|\pi_1| + |\pi_2| \geq 2n$.  Thus the only interesting cases
are to determine which pairs $(n,n)$, $(\overline n, \overline n)$, $(n, \overline n)$
are feasible.

The cominuscule recursion 
gives this answer to this question.
We use the computations~\eqref{E:hats}.
If $(n{-}1, \overline{n{-}1})$ is feasible for $Q^{2n-2}$, then
$|\overline n|_{n-1} + |n|_{\overline{n{-}1}} = 2 > 1$, and so 
$(n, \overline n)$
is infeasible for $Q^{2n}$, whereas $(n, n)$ and $(\overline n, \overline n)$ are 
feasible.  Similarly if $(n{-}1, n{-}1)$ and 
$(\overline{n{-}1}, \overline{n{-}1})$ are feasible for $L/Q$,
then
$|\overline n|_{n-1} + |\overline n|_{n-1} = 
|n|_{\overline{n{-}1}} + |n|_{\overline{n{-}1}} = 2 > 1$, and so 
$(n, n)$ and $(\overline n, \overline n)$ are infeasible, and $(n, \overline n)$
is feasible.  By induction, we see that if $n$ is odd,
$(n,n)$ and $(\overline n, \overline n)$ are feasible for $Q^{2n}$
and $(n, \overline n)$ is infeasible, and
vice-versa if $n$ is even.

\subsection{Type $B_n$, $G=SO_{2n+1}$, $G/P$ is an 
odd-dimensional quadric, $Q^{2n-1}$.}\label{S:odd-quadric}
The analysis of the odd-dimensional quadric is similar to the
even-dimensional quadric, in that $\frg/\frp$ is the defining
representation of $L=SO_{2n-1}$ and there is a single interesting
$L$-orbit on $\frg/\frp$ consisting of non-zero isotropic
vectors.  In the even-dimensional quadric, this
orbit gave the inequalities for determining feasibility in
the middle dimension.  For the odd-dimensional quadric, which
has no middle-dimensional cohomology, these inequalities are
redundant: they are all implied
by the basic codimension inequality $\sum_{i=1}^s |\pi_i| \leq 2n-1$.
Thus feasibility for $Q^{2n-1}$ is trivial, as
the only inequality needed is the basic codimension inequality.

\subsection{Type $C_n$, $G=Sp_{2n}$, $G/P$ is the Lagrangian Grassmannian.}
    \label{S:Lagrangian}
Suppose that $P$ is obtained by omitting the long root from the Dynkin diagram
for $C_n$.
Then $G/P=LG(n)$, the \Blue{{\em Lagrangian Grassmannian}} of isotropic $n$-planes
in $\CC^{2n}$, where $\CC^{2n}$ 
is equipped with a non-degenerate alternating bilinear form.
The Levi subgroup of $P$ is $GL_n$, and 
$\frg/\frp$ is 
the second symmetric power of the defining representation of
$GL_n$, that is, symmetric $n\times n$ matrices.
Its weights are $\{ e_i+e_j \mid 1\leq i\leq j\leq n\}$,
where $e_1,\dotsc,e_n$ are the standard orthonormal basis vectors of
$\CC^n=\frh^*$.

We identify $\Phi(\frg/\frp)$ with the cells of 
the staircase shape of height $n$.
Numbering the rows and columns in the standard way for matrices, 
the weight $e_i+e_j$ with $i\leq j$
corresponds to the cell in row $i$ and column $j$ in the staircase.
We write the coinversion set of a minimal coset representative 
$\pi\in W^P$ as a strict partition in the staircase, with the
inversion set its complement.
We use the strict partition of $\Invc(\pi)$ to represent elements
$\pi\in W^P$.
We display this for $n=7$ and 
$\pi=7521$, 
shading the inversions of $\pi$.
 \[
  \begin{picture}(180,80)(-3,0)
   \put(42,0){\includegraphics{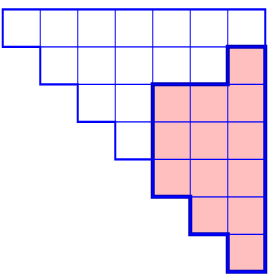}}
   \put( -3  ,59){$\Invc(\pi)$}  \put( 36,61){\vector(1,0){25}}
   \put(139.5,26){$\Inv(\pi)$} \put(136.5,29){\vector(-1,0){25}}
  \end{picture}
\]
The lowest root 
in $\Phi(\frg/\frp)$
is in the last row and the highest root is in the
first column.

The long roots 
in $\Phi(\frg/\frp)$
are $2e_1, \dotsc,2e_n$, which are
pairwise orthogonal.
Set $\alpha_i:=2e_{n+1-i}$,
 which is the $i$th root along the diagonal edge 
 of $\Phi(\frg/\frp)$ from the lower right. 
Then the vector $v$ has the form
$v_{\alpha_1}+ v_{\alpha_2} + \dotsb + v_{\alpha_r}$, 
The weights of $\frg\frl_n\cdot v$ are 
$\{e_i+e_j\mid n-r< j\}$,
and the subgroup $Q \in M(P)$ of $L=GL_n$ which is the stabilizer of 
$\frg\frl_n\cdot v$ is the stabilizer of
the $r$-dimensional linear subspace 
spanned by the last $r$ basis vectors $e_{n+1-r},\dotsc,e_n$.  
Thus $L/Q$ is the classical Grassmannian, $\Gr(r,n)$.
In this way, the weights of $\frg\frl_n\cdot v$ are the 
last $r$ columns of the
staircase and the weights of $\frz$ are the 
first $n{-}r$ columns and $\dim \frz=\binom{n{-}r{+}1}{2}$.
We 
show
this for $n=7$ and $r=3$.
 \[
  \begin{picture}(180,80)(10,0)
   \put(41,0){\includegraphics{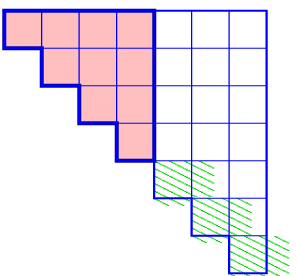}}
   \put(138.5,36){$\Phi(\frg\frl_n\cdot v)$} \put(135.5,39){\vector(-1,0){25}}
   \put( 10  ,58){$\Phi(\frz)$} \put(35,61){\vector(1,0){25}}
   \put( 49,26){$\alpha_{3}$} \put( 62,28){\vector(1,0){30}}
   \put( 60,15){$\alpha_{2}$} \put( 73,17){\vector(1,0){30}}
   \put( 71, 4){$\alpha_{1}$} \put( 84, 6){\vector(1,0){30}}
  \end{picture}
\]

\begin{figure}[htb]
 \[
  \begin{picture}(200,98)(0,-18)
   \put(8,0){\includegraphics{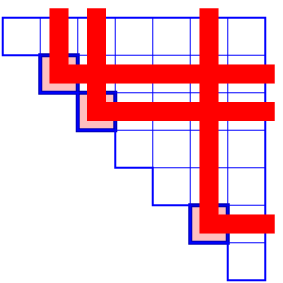}}
   \put( 1,58){1}   \put(11.5,47){2}   \put(22,36){3}
   \put(32.5,25){4}   \put(43,14){5}   \put(53.5, 3){6}
   \put(64,-8){7}

   \put( 32,-18){(a)}
   \put(136,-18){(b)}

   \put(112,0){\includegraphics{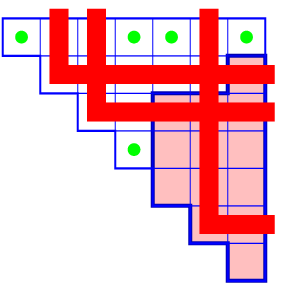}}
  \end{picture}
\]
\caption{$|\pi|_\lambda=5$ for $\pi=7521$ 
 and $\lambda=236$.\label{F:C}}
\end{figure}

Elements $\lambda\in W^Q$ act on $\Phi(\frg/\frp)$ by simultaneously shuffling 
rows and columns numbered 
$1,\dotsc,n{-}r$ with those numbered $n{+}1{-}r,\dotsc,n$.
This is equivalent to selecting $r$ boxes on the diagonal
corresponding to the images of the 
roots
$\alpha_1,\dotsc,\alpha_r$.
Then $\lambda\Phi(\frz)$ consists of weights which are orthogonal to each of the
selected weights, and are obtained by crossing out the row and column of
each selected box. 
This is displayed in Figure~\ref{F:C}(a) for $n=7$, $r=3$, and
when $\lambda$ selects the boxes in positions
$2$, $3$, and $6$.
After crossing out the rows and columns, 
$|\pi|_\lambda$ is the number of boxes in $\Invc(\pi)$ which are not crossed out.
We display this in Figure~\ref{F:C}(b) for $\pi=7521$ with the same
numbers $n$, $r$, and $\lambda$ as before.

We associate a minimal coset representative $\lambda\in W^Q$ for
$\Gr(r,n)$ to a selection of boxes on the diagonal 
in the same way as for columns
in Remark~\ref{Rem:rem}.
In our example, the selection of positions $2$, $3$, and $6$ gives the
inversion diagram for $\Gr(3,7)$.
\[
  \begin{picture}(63,47)
   \put(0,0){\includegraphics{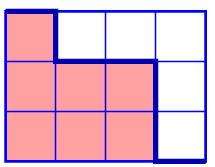}}
    \put(19,34){6}  \put(47,20){3} \put(47, 6){2} 
  \end{picture}
\]

\subsection{Type $D_{n+1}$, $G=SO_{2n+2}$, $G/P$ is the orthogonal Grassmannian,
  $OG(n{+}1)$.}
\label{S:Dn-spinor}
Suppose that $P$ is obtained by omitting one of the roots in the fork of the
Dynkin diagram for $D_{n+1}$.
Then $G/P$ is the orthogonal Grassmannian $OG(n{+}1)$ of isotropic
$n{+}1$-planes in $\CC^{2n+2}$, where $\CC^{2n+2}$ 
is equipped with a non-degenerate symmetric bilinear form.
The Levi subgroup of $P$ is $GL_{n+1}$, and $\frg/\frp$ is 
the second exterior power of the defining representation of $GL_{n+1}$,
that is, anti-symmetric $(n{+}1)\times(n{+}1)$-matrices.
Its weights are $\{e_i+e_j \mid 1\leq i< j\leq n+1\}$.

We identify $\Phi(\frg/\frp)$ with the 
cells of the staircase shape of height $n$.
Minimal coset representatives $\pi\in W^P$ are strict partitions
corresponding to $\Invc(\pi)$.
This is exactly the same as for the Lagrangian Grassmannian $LG(n)$;
not only do these two cominuscule flag varieties have
Schubert positions indexed by the same set (of strict partitions), but a
Schubert position $\pi_1,\dotsc,\pi_s$ is feasible for $LG(n)$ if and
only if $\pi_1,\dotsc,\pi_s$ is feasible for 
$OG(n{+}1)$. 
Despite this similarity, the minuscule recursion is different for $LG(n)$ and
for $OG(n{+}1)$. 

Numbering the rows of the staircase from 1 to $n$ with 1 the longest 
row, and the columns $2$ to $n{+}1$ with $n{+}1$ the longest column, 
the weight $e_i+e_j$ with $i< j$ 
corresponds to the cell in row $i$ and column $j$ in the staircase. 
Every root 
in $\Phi(\frg/\frp)$ is long.
Set 
%
%
$\alpha_i:=e_{n+2-2i}+e_{n+3-2i}$, 
 which is the 
$(2i{-}1)$st
root along the diagonal edge of 
  $\Phi(\frg/\frp)$ from the lower right.
Then the vector $v$ has the form
\[
   v_{\alpha_{1}}+ v_{\alpha_{2}} + \dotsb + v_{\alpha_{r}}\ .
\]
The weights of $\frg\frl_{n+1}\cdot v$ are $\{e_i+e_j\mid n+1-2r<j\}$,
and the subgroup $Q$ of $L=GL_{n+1}$ which stabilizes $\frg\frl_{n+1}\cdot v$ 
is the subgroup stabilizing the $2r$-dimensional linear subspace 
spanned by the last $2r$ basis vectors, $e_{n+2-2r},\dotsc,e_{n+1}$.
Thus 
$L/Q$ is an ordinary Grassmannian $\Gr(2r,n)$ of even-dimensional subspaces.
In this way, the weights of $\frg\frl_{n+1}\cdot v$ are the last $2r$ columns of the
staircase and the weights of $\frz$ are the first $n{-}2r$ columns.
We show this for $n=8$ and $r=2$.
 \[
 \begin{picture}(200,90)(10,0)
   \put(41,0){\includegraphics{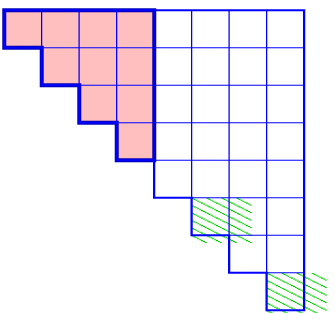}}
   \put(149.5,46){$\Phi(\frg\frl_{n+1}\cdot v)$} \put(146.5,49){\vector(-1,0){25}}
   \put( 10,68){$\Phi(\frz)$} \put(35,71){\vector(1,0){25}}
   \put( 60,26){$\alpha_{2}$}             \put( 73,28){\vector(1,0){30}}
   \put( 82,4){$\alpha_{1}$}             \put( 95,6){\vector(1,0){30}}
  \end{picture}
\]

\begin{figure}[htb]
 \[
  \begin{picture}(230,119)(-6,-18)
   \put(10,2){\includegraphics{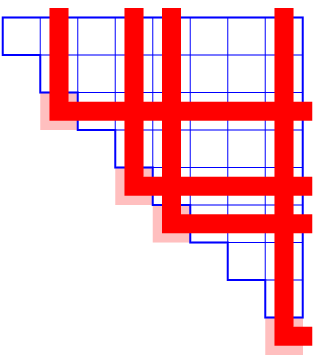}}
   \put(-6,81){1}   \put( 5,70){2}   \put(16,59){3}
   \put(27,48){4}   \put(38,37){5}   \put(49,26){6}
   \put(60,15){7}   \put(71, 4){8}   \put(82,-7){9}

   \put( 39,-18){(a)}
   \put(163,-18){(b)}

   \put(132,5){\includegraphics{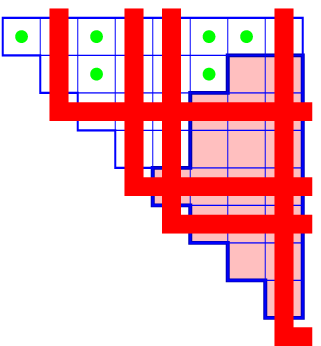}}
  \end{picture}
\]
\caption{$|\pi|_\lambda=6$ for $\pi=8532$ and $\lambda=3569$.\label{F:D}}
\end{figure}

Elements $\lambda\in W^Q$ act on $\Phi(\frg/\frp)$ by permuting the indices of
the weights $e_i+e_j$.
Since 
\[
   (e_i+e_j, e_k+e_l)\ =\ |\{i,j\}\cap\{k,l\}|\,,
\]
we obtain the weights of $\lambda\Phi(\frz)$ as follows.
The diagonal positions in row and column $i$ for $i=1,\dotsc,n{+}1$
lie outside the staircase.
Then $\lambda$ selects $2r$ of these positions, 
and as before, we cross out the rows and columns 
of these $2r$ positions.
This is displayed in Figure~\ref{F:D}(a) for $n=8$, $r=2$, and $\lambda=3569$.
Then $|\pi|_\lambda$ counts the boxes in $\Invc(\pi)$ which are not crossed out.
We display this in Figure~\ref{F:D} for $\pi=8532$ with the same numbers
$n$, $r$, and $\lambda$ as before.
For this case, $|\pi|_\lambda=6$.
We associate a minimal coset representative $\lambda\in W^Q$ for
$\Gr(2r,n{+}1)$ to a selection of boxes on the diagonal 
in the same way as for columns in Remark~\ref{Rem:rem}.

We note that the inequalities for $OG(n{+}1)$ are quite different than the inequalities of 
Section~\ref{S:Lagrangian} for the Lagrangian Grassmannian $LG(n)$, despite their having
the same sets of solutions.

\subsection{Type $E_6$, $G/P$ is the Cayley plane $\OP^2$.}
\label{S:E6}

\begin{figure}[htb]
\[
 \begin{picture}(260,165)(0,5)
  \put(  0,0){\includegraphics{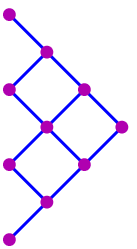}}
   \put( 25,5){$\mathcal{E}_4$}
  \put( 80,0){\includegraphics{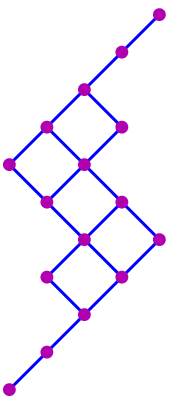}}
   \put(110,5){$\mathcal{E}_5$}
  \put(180,0){\includegraphics{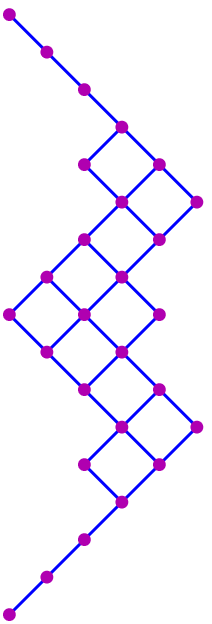}}
   \put(215,5){$\mathcal{E}_6$}
 \end{picture}
\]
\caption{Lattices for the exceptional cominuscule flag varieties.\label{F:exceptional}}
\end{figure}

This 
is in many ways similar to the even-dimensional quadric.
Here, the parabolic subgroup $P$ is obtained by omitting the rightmost node of the 
Dynkin diagram of $E_6$, as shown in Table~\ref{Ta:One}.
Its Levi subgroup $L$ has semisimple part 
$\text{Spin}_{10}$ (type $D_5$), and the flag 
variety $G/P$ is the
even-dimensional Cayley plane $\OP^2$.
The lattice $\Phi(\frg/\frp)$ is the poset $\mathcal{E}_5$ of
Figure~\ref{F:exceptional}.
Thus $W^P$ is the set of (lower) order ideals in $\mathcal{E}_5$, where
$\pi \in W^P$ corresponds to the order ideal $\Inv(\pi)$.

%
%
%

The tangent space $\frg/\frp$ is the $16$-dimensional spinor representation
of $L$.  
As in Section~\ref{S:even-quadric}, $M(P)$ consists of a single
parabolic subgroup $Q$, where $L/Q = OG(5)$. 
The $H$-fixed points
on $L/Q$ are the images of the weight spaces of $\frg/\frp$, and
thus $W^Q$ is canonically identified with $\mathcal{E}_5$.

If $\alpha$ is the simple root defining $P$, then $\Phi(\frz)$ is
the orthogonal complement $\alpha^\perp$ to $\alpha$ in $\Phi(\frg/\frp)$,
which consists of $5$ roots.  Moreover for $\lambda \in W^Q$, we have
$\lambda \Phi(\frz) = \lambda^\perp$ is the orthogonal complement to
$\lambda$ in $\Phi(\frg/\frp)$.  Consequently, viewing $\pi$ as
an order ideal in $\mathcal{E}_5$, and $\lambda$ as an element of 
$\mathcal{E}_5$,
we have the following formula:
\[
  |\pi|_\lambda 
  \ =\  \big|\{\beta \in \Phi(\frg/\frp)
    \mid \beta \notin \pi,\ \beta \perp \lambda\}\big|\,,
\]
and the inequalities~\eqref{E:w-ineq} are
\[
  \sum_{i=1}^s |\pi_i|_{\lambda_i}\ \leq\ 5\,.
\]

Note that the weight lattice $\Phi(\frl/\frq)$ is isomorphic
to $\mathcal{E}_4$.
There is a unique isomorphism from $\mathcal{E}_5$ to 
$J(\mathcal{E}_4)$.  Thus we can view each $\lambda \in W^Q$
as an order ideal in $\mathcal{E}_4$, which is 
a strict partition inside a staircase diagram.  This allows us
to continue the recursion with $OG(5)$, as discussed in 
Section~\ref{S:Dn-spinor}.

\subsection{Type $E_7$, $G/P$ is $G_\omega(\OO^3,\OO^6)$.}
  
The parabolic subgroup $P$ is obtained by omitting the rightmost node of the 
Dynkin diagram of $E_7$, as shown in Table~\ref{Ta:One},
Its Levi subgroup $L$ has type $E_6$ and the flag 
variety $G/P = G_\omega(\OO^3,\OO^6)$.
The lattice $\Phi(\frg/\frp)$ is the poset $\mathcal{E}_6$, so
that $\pi \in W^P$ corresponds to an order ideal in $\mathcal{E}_6$,
via its inversion set.

The tangent space $\frg/\frp$ is the $27$-dimensional minuscule
representation of $E_6$.  This has two interesting orbits.  The smallest
is the orbit through $v = v_\alpha \in \frg/\frp$, where $\alpha$ is 
the simple root defining $P$.  It is $17$-dimensional, and gives rise to 
the parabolic subgroup $Q \subset L$ which obtained by omitting the rightmost node 
of the $E_6$ Dynkin diagram.  The second orbit is $26$-dimensional,
and is the orbit through 
$v = v_\alpha + v_{\alpha_2}$, where $\alpha_2 \in \Phi(\frg/\frp)$ is 
the (unique) lowest root orthogonal to $\alpha$.  This orbit gives
rise to the parabolic subgroup $Q \subset L$ obtained by omitting the leftmost
node of the $E_6$ Dynkin diagram.  Thus in both cases $L/Q$ is isomorphic 
to the Cayley plane $\OO\PP^2$, but these two manifestations of the Cayley 
plane give rise to different inequalities.  
(This also occurs for $LG(n)$, 
where we have different inequalities coming from
isomorphic varieties $\Gr(r,n)$ and $\Gr(n{-}r, n)$.)

As in Section~\ref{S:E6}, the Schubert positions for $\OO\PP^2$
correspond to order ideals in 
$\mathcal{E}_5$.  Since $J(\mathcal{E}_5)$ is canonically isomorphic
to $\mathcal{E}_6$, we will now 
identify $W^Q$ with $\mathcal{E}_6$.

For the smaller orbit, $\Phi(\frz)$ is the orthogonal complement 
$\alpha^\perp$ to $\alpha$ in $\Phi(\frg/\frp)$, which consists of $10$ roots.  Thus
viewing $\pi$ as
an order ideal in $\mathcal{E}_6$, and $\lambda$ as an element of 
$\mathcal{E}_6$,
we have the following formula:
\[
  |\pi|_\lambda 
  = \big|\{\beta \in \Phi(\frg/\frp)
    \mid \beta \notin \pi,\ \beta \perp \lambda\}\big|\,.
\]
and the inequalities~\eqref{E:w-ineq} for this orbit are
\[
  \sum_{i=1}^s |\pi_i|_{\lambda_i}\ \leq\ 10\,.
\]

For the larger orbit, $\Phi(\frz)$ is the orthogonal complement to
$\{\alpha, \alpha_2\}$, which consists of highest root 
in $\Phi(\frg/\frp)$.
Let $\lambda \mapsto \widehat{\lambda}$ denote the unique order reversing
involution on $\mathcal{E}_6$.  Then $\Phi(\frz) = \widehat{\alpha}$, and 
in general $\lambda\Phi(\frz)$ is the single root $\widehat{\lambda}$.
Thus we have
\[
  |\pi|_\lambda 
 \  =\  
   \begin{cases}
   0 &\text{if $\widehat{\lambda} \in \pi$} \\
   1 &\text{otherwise}
   \end{cases}
\]
and the inequalities~\eqref{E:w-ineq} for this orbit are
\[
  \sum_{i=1}^s |\pi_i|_{\lambda_i}\ \leq\ 1\,.
\]

\section{Comparison with other inequalities}\label{S:five}

We first 
discuss how the classical Horn inequalities arise from the
inequalities of Theorem~\ref{T:necessary} and how to modify the proof of
Theorem~\ref{Theorem} to prove their sufficiency.
Next, we show how to use Proposition~\ref{P:FeasibleSP} to derive a
different set of necessary inequalities for feasibility on $G/P$, which
we call the naive inequalities.
When $G/P$ is the classical Grassmannian, these include the Horn
inequalities and were essentially derived by
Fulton~\cite[Section~1]{Fu00a}. 

Our derivation of naive inequalities generalizes Theorem~36 of Belkale
and Kumar in~\cite{BK04}.
While their subset is a proper subset of the
inequalities~\eqref{E:generalineq} from 
Theorem~\ref{T:necessary}, it includes none of the sufficient
inequalities~\eqref{E:w-ineq}. 

Finally, we explain these naive inequalities in detail for the
Lagrangian Grassmannian, which shows they are quite different than the
inequalities of Theorem~\ref{Theorem}, as given in
Section~\ref{S:Lagrangian}. 
We conjecture that the naive inequalities are sufficient to determine
feasibility for the Lagrangian Grassmannian. 
We have verified this conjecture for $s=3$ and $n\leq 8$.

\subsection{Horn inequalities}\label{S:Horn}

Schubert classes $\sigma_\mu$ in the cohomology of the
Grassmannian $\Gr(k,n)$ are traditionally indexed by partitions $\mu$,
which are weakly decreasing sequences of non-negative integers
\[
   \mu\ \colon\ n{-}k\ \geq\ \mu^1\ \geq\ \mu^2\ 
               \geq\ \dotsb\ \geq\ \mu^k\ \geq\ 0\,.
\]
Write $|\mu|$ for the sum $\mu^1+\dotsb+\mu^k$.
The partition $\mu$ associated to a Schubert position $\pi$ is
essentially its coinversion set $\Invc(\pi)$.
Specifically, $\mu^i$ is the number of positive roots of the form
$e_j-e_i$ which are coinversions.
With the conventions of Section~\ref{S:An}, the Ferrers diagram of
$\mu$ is the reflection of $\Invc(\pi)$ across a vertical line.

Let $\mu^t$ denote the conjugate partition to $\mu$, whose Ferrers
diagram is obtained by transposing the Ferrers diagram of $\mu$.
Note that if $\mu$ indexes a Schubert class for $\Gr(k,n)$,
then $\mu^t$ indexes a Schubert class for $\Gr(n{-}k,n)$.

Given Schubert positions $\mu_1,\dotsc,\mu_m$ and $\nu$ for $\Gr(k,n)$,
we say that $\sigma_\nu$ \Blue{{\em occurs} in $\prod_{i=1}^m \sigma_{\mu_i}$}
if, when we expand the product in the basis of Schubert
classes, $\sigma_\nu$ occurs with a non-zero coefficient.
Necessarily, we must have the codimension condition
\[
  |\nu|\ =\ |\mu_1|+|\mu_2|+\dotsb+|\mu_m|\,.
\]

If $\mu$ is a partition indexing a Schubert position for $\Gr(k,n)$,
and $\kappa \colon k-r \geq \kappa^1 \geq \dotsc \geq \kappa^r \geq 0$
is a partition for $\Gr(r,k)$, let
\[
  \kappa[a]\ :=\ a+\kappa^{r+1-a}\qquad \text{and}
\qquad
  |\mu|^\kappa\ :=\ \sum_{a=1}^r \mu^{\kappa[a]}\,.
\]
We recall the Horn recursion for $\Gr(k,n)$,
following Fulton~\cite[Theorem~17(1)]{Fu00}. 

\begin{prop}\label{P:Classical_Horn}
  Let $\mu_1,\dotsc,\mu_m$ and $\nu$ be Schubert positions for $\Gr(k,n)$
  with $|\nu| = |\mu_1|+\dotsb+|\mu_m|$.
  The following are equivalent.
\begin{enumerate}
 \item[(i)] $\sigma_\nu$ occurs in $\prod_{i=1}^m \sigma_{\mu_i}$.

 \item[(ii)] The inequality
 \begin{equation}\label{Eq:C_Horn}
   \sum_{i=1}^m |\mu_i|^{\kappa_i}\ \geq\ 
                |\nu|^\theta
 \end{equation}
 holds for all Schubert positions 
 $\kappa_1,\dotsc,\kappa_m$ and $\theta$ for $\Gr(r,k)$ such that 
 $\sigma_\theta$ occurs in $\prod_{i=1}^m \sigma_{\kappa_i}$, and all
 $1\leq r<k$.
\end{enumerate}
\end{prop}

The proof of Theorem~\ref{Theorem} can be modified to
prove Proposition~\ref{P:Classical_Horn}.

As we saw in Section~\ref{S:An}, the  semisimple part of the Levi subgroup 
is a product $L^{ss} = L_0 \times L_1$.  
Rather than study the tangent space $\frl \cdot v$ to the $L$-orbit
through $v$, we 
instead study the tangent space $\frl_1 \cdot v$ to
the $L_1$-orbit through $v$.
Then Lemma~\ref{L:Inductive_engine} is 
true under this substitution for the following reason.
Let $\phi_1$ denote the
new quotient map $\phi_1\colon(\frg/\frp) \to (\frg/\frp)/(\frl_1 \cdot v)$.
We know from Lemma~\ref{L:Inductive_engine} (as originally stated) that 
the intersection 
$\bigcap_{i=1}^s l_iT_{\pi_i}$ is transverse if and only if
$\bigcap_{i=1}^s \phi(l_iT_{\pi_i})$ is transverse. 
But since $\phi$ factors through $\phi_1$, by
Proposition~\ref{P:Transversality}(i) these are transverse if and only
$\bigcap_{i=1}^s \phi_1(l_iT_{\pi_i})$ is transverse.  
The rest of
the proof proceeds very much as written (although most of the Appendix
is unnecessary since this is type $A$).  
We deduce that by using only
one factor of $L^{ss}$, one obtains a set of necessary and sufficient 
inequalities for feasibility on Grassmannians, different from those
of Theorem~\ref{Theorem}.

These inequalities turn out to be the classical
Horn inequalities.  
To see this, we adopt some of the notation of Section~\ref{S:An},
identifying $\frg/\frp$ with $\Hom(\CC^k,\CC^{n-k})$ and 
$L^{ss}$ with $SL_k(\CC)\times SL_{n-k}(\CC)$, where $L_1=SL_{n-k}(\CC)$.
If $v$ has the form~\eqref{Eq:vee}, then 
\[
  \Phi(\frs\frl_{n-k}\cdot v)\ =\ 
  \Phi(\frl_1\cdot v)\ =\ 
  \{ e_j{-}e_i\mid k-r< i\leq k< j\leq n\}\,.
\]
Thus $\Phi(\frz)$ is the upper $(k{-}r)\times(n{-}k)$ rectangle in
$\Phi(\frg/\frp)$, so that $\dim \frz = (k{-}r)(n{-}k)$.
We display this when $n=11$, $k=5$, and $r=2$.
\[
  \begin{picture}(80,58)
   \put(0,0){\includegraphics{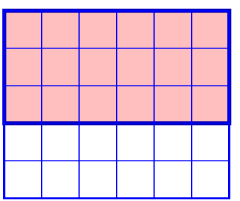}}
   \put(86.5,36){$\Phi(\frz)$} \put(83.5,39){\vector(-1,0){25}}
   \put(86.5, 4){$\Phi(\frs\frl_{n-k}\cdot v)$} 
    \put(83.5, 7){\vector(-1,0){25}}
  \end{picture}
\]

The subgroup $Q$ which is the stabilizer of $\frs\frl_{n-k}\cdot v$ is
obtained by further omitting the node at position $k{-}r$ in the Dynkin
diagram for $L$. 
Thus $L/Q$ is isomorphic to $\Gr(k{-}r,k)$.
Elements $\lambda\in W^Q$ act on $\Phi(\frg/\frp)$ by shuffling the $r$
rows which do not meet $\Phi(\frz)$ with those that do.
As before, $|\pi|_\lambda$ is the number of boxes in
$\Invc(\pi)$ which remain after crossing out the images of the rows not
in $\frz$.
For example, when  $n=11$, $k=5$, $\pi=1367\,\ten\,24589\,\eleven$, $r=2$, and 
we select rows 2 and 4 from the top, we see that  $|\pi|_\lambda=10$.
\[
    \includegraphics{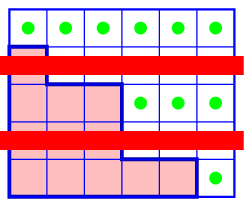}
\]

The preceding discussion shows that we have the following
recursion for top-degree Schubert positions (the analog of
Corollary~\ref{C:topdegree}). 

\begin{prop}[Horn recursion]\label{P:new_Horn}
 Let $\pi_1,\dotsc,\pi_s$ be a top-degree Schubert position for
 $\Gr(k,n)$.
 Then $\pi_1,\dotsc,\pi_s$ is feasible if and only if for every
 $1\leq r<k$ and every feasible top-degree Schubert position 
 $\lambda_1,\dotsc,\lambda_s$ for $\Gr(k{-}r,k)$, we have
 \begin{equation}\label{Eq:new_Horn} 
  \sum_{i=1}^s |\pi_i|_{\lambda_i}\ \leq\ (k-r)(n-k)
 \end{equation}
\end{prop}

Finally, we show that
the two recursions in
Propositions~\ref{P:Classical_Horn} and~\ref{P:new_Horn} are identical. 

Let $\mu$ be the partition associated to $\pi$ and $\kappa^t$ be the partition
associated to $\lambda$; thus $\kappa$ is a partition for $\Gr(r,k)$.
If we compare the definition of $|\pi|_\lambda$ with
Remark~\ref{Rem:rem}, which explains how to associate an 
inversion diagram to the rows selected by $\lambda\in W^Q$, 
we see that 
 \begin{equation}\label{Eq:relate}
   |\pi|_\lambda\ =\ |\mu|\ -\ |\mu|^\kappa\,.
 \end{equation}

Let $\nu$ be a partition for a Schubert position for
$\Gr(k,n)$.
The dual partition $\widehat{\nu}$ defined by
\[
    \widehat{\nu}\,^a\ =\ n-k\ -\ \nu^{k+1-a}\,,
\]
has the property that $|\nu|+|\widehat{\nu}|=k(n{-}k)$ and 
\[
   \int_{\Gr(k,n)} \sigma_{\mu}\sigma_{\widehat{\nu}}
   \ =\ \left\{\begin{array}{rcl}1&\ &\mbox{if } \mu=\nu\\
                                 0&\ &\mbox{otherwise}\,.
         \end{array}\right.
\]
Thus $\sigma_\nu$ appears in $\prod_{i=1}^m\sigma_{\mu^i}$ if and only if
$\mu^1,\dotsc,\mu^m, \widehat{\nu}$ is a feasible top-degree Schubert
position.

The reader can easily verify that 
$|\nu|^\theta = r(n-k) - |\widehat{\nu}|^{\widehat{\theta}}$.  
Thus~\eqref{Eq:C_Horn} becomes
 \begin{equation}\label{Eq:penultimate}
   \sum_{i=1}^m |\mu_i|^{\kappa_i}
   \ +\ |\widehat{\nu}|^{\widehat{\theta}} \geq\ 
   r(n{-}k)\,.
 \end{equation}
Since $\mu_1,\dotsc,\mu_m,\,\widehat{\nu}$ is a top-degree Schubert
position for $\Gr(k,n)$,
\[
   |\mu_1|+|\mu_2|+\dotsb+|\mu_m|+|\widehat{\nu}|\ =\ k(n{-}k)\,.
\]
We subtract~\eqref{Eq:penultimate} from this, setting $s:=m+1$, 
$\mu^s:=\widehat{\nu}$, and $\kappa^s:=\widehat{\theta}$, to obtain
\[
   \sum_{i=1}^s 
     \Bigl(|\mu_i|\ -\ |\mu_i|^{\kappa_i}\Bigr)
      \ \leq\ (k{-}r)(n{-}k)\,.
\]
If the partition $\mu_i$ corresponds to the representative $\pi_i\in W^P$
and the partition $\kappa_i^t$ to the representative $\lambda_i\in W^Q$,
then, by~\eqref{Eq:relate}, this is just the condition~\eqref{Eq:new_Horn}.

\subsection{Naive inequalities}

 Recall the situation of Proposition~\ref{P:FeasibleSP}.
 We have parabolic subgroups $R\subset P\subset G$ and a correspondence
 between Schubert positions $\lambda$ for $P/R (=L/Q)$,
 $\pi$ for $G/P$, and $\lambda\pi$ for $G/R$.

 Suppose that $P'$ is another parabolic subgroup of $G$ which contains
 $R$.
 The image of the Schubert variety $X_{\lambda\pi}$ of $G/R$ under the 
 projection to $G/P'$ is a (translate of a) Schubert variety 
 $X_{\pi'}$ of $G/P'$.
 Write $\|\pi\|_{\lambda}$ for the codimension of $X_{\pi'}$ in $G/P'$.
 We intentionally suppress the dependence of $\pi'$ on $\lambda$ and
 of $\|\pi\|_{\lambda}$ on $P'$.
 We use Proposition~\ref{P:FeasibleSP}, which relates feasibility for Schubert
 problems on different flag varieties, to obtain necessary inequalities which
 hold for feasible Schubert problems on $G/P$.

\begin{thm}\label{T:naive}
 Suppose that $\pi_1,\dotsc,\pi_s$ is a feasible Schubert position for
 $G/P$.
 Given parabolic subgroups $R\subset P'$ of $G$ with 
 $R\subset P$ and any feasible Schubert position
 $\lambda_1,\dotsc,\lambda_s$ for $P/R$, the Schubert position
 $\pi'_1,\dotsc,\pi'_s$ is feasible for $G/P'$.
 In particular, any necessary inequalities for feasibility on
 $G/P'$ give inequalities on the original Schubert position
 $\pi_1,\dotsc,\pi_s$ for $G/P$.
 For example, the basic codimension inequality for $\pi'_1, \dotsc, \pi'_s$
gives
 \begin{equation}\label{Eq:Necess_codim}
    \sum_i \|\pi_i\|_{\lambda_i}\ \leq\ \dim G/P'\,.
 \end{equation}
\end{thm}

\begin{proof}
 By Proposition~\ref{P:FeasibleSP}(ii),
 $\lambda_1\pi_1,\dotsc\lambda_s\pi_s$ is a feasible Schubert position
 for $G/R$, and so by Proposition~\ref{P:FeasibleSP}(i), 
 $\pi'_1,\dotsc,\pi'_s$ is a feasible Schubert position for $G/P'$.
 The rest is immediate.
\end{proof}

\begin{remark}\label{R:Bel-Kum}
Belkale and Kumar~\cite[Theorem 36]{BK04} use similar ideas to
also derive~\eqref{Eq:Necess_codim}. 
When $P'\cap P = R$, they express these inequalities in a form
similar to~\eqref{Eq:New_Horn}, in terms of counting 
roots~\cite[inequality (58)]{BK04}.
In fact, these are the inequalities of Theorem~\ref{T:necessary}, when
$\frs=\frn_{P'}\cap\frn_P$. 
Since this is almost never the center
of $\frn_R$, none of the inequalities of Belkale and Kumar have
the form~\eqref{E:w-ineq}.

We note that the inequalities~\eqref{Eq:Necess_codim} are 
always a subset of the
necessary inequalities of Theorem~\ref{T:necessary}.
The verification of this assertion is left as an exercise.

\end{remark}

\begin{remark}\label{R:Plethora}
 Theorem~\ref{T:naive} gives a method to generate many necessary
 inequalities for feasibility on different flag varieties $G/P$.
 For example, in type $A$ we can take $R=P$ and let $P'$ be any maximal
 parabolic subgroup containing $P$.
 Then $G/P'$ is a classical Grassmannian and Theorem~\ref{T:naive} shows
 how to pull back the Horn inequalities for $G/P'$ to obtain
 inequalities for $G/P$. 

 If $G$ has type $A$, $B$, $C$, or $D$ and $P$ is a
 maximal parabolic subgroup, then we can select $R\subset P$ so that $P/R$
 is a classical Grassmannian.
 If $P'$ is a different parabolic subgroup of $G$ which contains $R$, then
 the codimension inequalities~\eqref{Eq:Necess_codim} give necessary
 inequalities for feasibility on $G/P$ which are indexed by feasible
 Schubert problems on a classical Grassmannian $P/R$.

 We invite the reader to check that in type $A$, this last procedure is
 yet another method for deriving the necessity of the Horn inequalities.
 In fact, Fulton essentially did just that in~\cite[Section~1]{Fu00a}.

 We also invite the reader to use Theorem~\ref{T:naive} to generate
 even more necessary inequalities for feasibility on flag varieties $G/P$.
 We believe that it is an interesting and worthwhile project to investigate
 these naive Horn-type inequalities on other flag varieties.
 For example, for which flag varieties is (a natural subset of) the set of
 all such naive inequalities sufficient to determine feasibility?
 In the next section, we examine a subset of these in detail for the
 Lagrangian Grassmannian, showing that they are in general different than
 the necessary and sufficient inequalities derived in Section~\ref{S:Lagrangian}.
\end{remark}

\subsection{Naive inequalities for the Lagrangian Grassmannian}
We express codimension inequalities~\eqref{Eq:Necess_codim} of
Theorem~\ref{T:naive} for the Lagrangian Grassmannian in a form similar to the
inequalities of Corollary~\ref{C:topdegree}. 

\begin{thm}\label{Naive_LG}
  Let $\pi_1,\dotsc,\pi_s$ be a feasible top-degree Schubert position for the
  Lagrangian Grassmannian $LG(n)$.
  Then, for any feasible Schubert positions $\lambda_1,\dotsc,\lambda_s$ for
  $\Gr(r,n)$, we have
 \begin{equation}\label{Eq:LG_Naive}
   \sum_{i=1}^s |\pi_i|_{\widehat{\lambda}_i^t}\ \geq\ \binom{r{+}1}{2}\ .
 \end{equation}
\end{thm}

Here, $\widehat{\lambda}^t$ is the conjugate of the
dual Schubert position to $\lambda$, as in
Section~\ref{S:Horn}. 

Note that $\lambda_1, \dotsc, \lambda_s$ is feasible for $\Gr(r,n)$ if
and only if $\lambda_1^t, \dotsc, \lambda_s^t$ is feasible for $\Gr(n{-}r,n)$.
Thus the inequalities~\eqref{Eq:LG_Naive} may be rewritten
\[
   \sum_{i=1}^s |\pi_i|_{\widehat{\lambda}_i}\ \geq\ \binom{n{-}r{+}1}{2}\ .
\]
These bear a striking similarity to the inequalities of
Corollary~\ref{C:topdegree}
for the Lagrangian Grassmannian, which
by the discussion in Section~\ref{S:Lagrangian}, have the form
\[
   \sum_{i=1}^s |\pi_i|_{\lambda_i}\ \leq\ \binom{n{-}r{+}1}{2}\ ,
\]
and are indexed by the same set as the necessary inequalities of
Theorem~\ref{Naive_LG}.  In fact these inequalities 
are quite different.
Not only does the inequality go in the opposite direction, but the 
terms
$|\pi_i|_{\widehat{\lambda}_i}$ and $|\pi_i|_{\lambda_i}$ are 
unrelated quantities.

\begin{proof}
 Let $G=Sp_{2n}(\CC)$ and $P_k$ be the maximal parabolic subgroup corresponding
 to the $k$th simple root from the right
 end of the Dynkin diagram of $C_n$ as
 shown in Table~\ref{Ta:One}.
 Then $G/P_k$ is a space of isotropic $k$-dimensional linear subspaces of
 a $\CC^{2n}$ which is equipped with a non-degenerate alternating 
 bilinear form and
 $\dim G/P_k=2k(n-k) + \binom{k{+}1}{2}$.

 We consider the codimension inequalities~\eqref{Eq:Necess_codim} of
 Theorem~\ref{T:naive} for $G/P_n=LG(n)$, the Lagrangian Grassmannian, 
 $P'=P_{n-r}$, and $R=P_n\cap P_{n-r}$.
 Let $\pi$ be a Schubert position for $LG(n)$ and $\lambda$ a Schubert position
 for $P_n/R=\Gr(r,n)$.  (Note: it is consistent with 
 the conventions established in Section~\ref{S:Lagrangian}, 
 to call this $\Gr(r,n)$, rather than $\Gr(n{-}r,n)$.)

 We will show
 \begin{equation}\label{Eq:critical}
   \| \pi \|_\lambda\ =\ |\pi|\ +\ |\lambda|\ -\ |\pi|_{\widehat{\lambda}^t}\,.
 \end{equation}
 Then~\eqref{Eq:LG_Naive} will follow, for 
 \begin{eqnarray*}
  \sum_{i=1}^s |\pi_i|_{\widehat{\lambda}_i^t}& =& 
     \sum_{i=1}^s  |\pi_i|\ +\ \sum_{i=1}^s |\lambda_i|\ -\ 
       \sum_{i=1}^s \|\pi_i\|_{\lambda_i}\\
      & \geq& \binom{n{+}1}{2} + r(n-r) - 2r(n-r) - \binom{n{-}r{+}1}{2}
     \ =\ \binom{r{+}1}{2}\,.
 \end{eqnarray*}
 Indeed, $\sum_i|\pi_i|=\binom{n{+}1}{2}$ and
 $\sum_i|\lambda_i|=r(n{-}r)$, as
 these are top-degree Schubert positions, and the inequality comes 
 from the negative of inequality~\eqref{Eq:Necess_codim}.

 We deduce~\eqref{Eq:critical} using a uniform combinatorial model for Schubert
 positions in these flag varieties, which may be found in~\cite{FP98}.
 Schubert positions $w$ for $G/P_k$ are represented by increasing sequences of
 integers
 \[
   w\ \colon\ 1\leq w^1 < w^2 <\dotsb< w^k\leq 2n\,,
 \]
 where we do not have $w^i+w^j=2n+1$ for any $i,j$.
 (The corresponding Schubert variety consists of those isotropic $k$-planes 
 $V$ where 
  $\dim V\cap F^{2n+1-w^i}\geq i$, 
 where $F^1,F^2,\dotsc$ is 
 a fixed isotropic flag with $i=\dim F^i$.)
 Then
 \begin{equation}\label{Eq:col_sum}
   |w|\ :=\ \left(\sum_{j=1}^k w^j-j \right)\ 
     -\ \bigl|\{ a<b\mid w^a+w^b>2n+1\}\bigr|\,.
 \end{equation}

Let $k=n$ and recall the conventions of Section~\ref{S:Lagrangian} 
for drawing coinversion sets for $\pi \in W^{P_n}$ as
strict partitions in the staircase shape with diagonal boxes
(and hence rows and columns) labeled $1,\dotsc,n$.
Let $w$ be the increasing sequence of integers corresponding to 
$\pi \in W^{P_n}$. The correspondence is such that
$\pi$ has a coinversion in position $(n{+}1{-}a,n{+}1{-}b)$ if and only
$w^a + w^b > 2n+1$.
The term $w^j - j$ in~\eqref{Eq:col_sum} counts the number
of coinversions in the hook through row and column $n{+}1{-}j$, while
$\bigl|\{ a<b \mid w^a+w^b>2n+1\}\bigr|$ is the total number of off-diagonal
coinversions, which are counted twice in the sum.

Let $\kappa^t$ be the partition corresponding to a Schubert position
$\lambda$ for $\Gr(r,n)$; thus $\kappa$ indexes a Schubert position for
$\Gr(n{-}r,n)$.   Recall that $\kappa[a] := a + \kappa^{n-r+1-a}$.  
If we lift a Schubert position $w$ for $G/P_n$ to
 $G/R$ using $\kappa^t$ and then project to $G/P_{n-r}$, we obtain the 
Schubert position
\[
   w'\ :=\ 
   w^{\kappa[1]}\ <\ w^{\kappa[2]}\ <\ \dotsb\ <\ w^{\kappa[n{-}r]}\,,
\]
and so
\[
   \|\pi\|_\lambda\ =\  |w'|\ =\ 
    \left(\sum_{j=1}^{n-r} w^{\kappa[j]} -j \right)\ -\ 
     \bigl|\{a<b\mid w^{\kappa[a]}+w^{\kappa[b]}>2n+1\}\bigr|\ .
\]

  Consider $\|\pi\|_\lambda-|\lambda| 
   = \|\pi\|_\lambda - \sum_{j=1}^s\kappa_{n-r+1-j}$, which is 
\[
   \left(\sum_{j=1}^{n-r} w^{\kappa[j]} -\kappa[j] \right)\ -\ 
     \bigl|\{a<b\mid w^{\kappa[a]}+w^{\kappa[b]}>2n+1\}\bigr|\,.
\]
 From the discussion interpreting the terms of~\eqref{Eq:col_sum} for $LG(n)$,
 it follows that this sum is the number of coinversions of $\pi$
 which lie in the hooks through rows and columns indexed
 $n{+}1{-}\kappa[j]$, for $j=1,\dotsc,n{-}r$.  The subtracted term
 $\bigl|\{a<b\mid w^{\kappa[a]}+w^{\kappa[b]}>2n+1\}\bigr|$ is the number
 of such coinversions counted twice by the sum.
 From the definition of $\widehat{\kappa}$ these are the rows and columns 
 indexed
 by $\widehat{\kappa}[j]$, for $j=1,\dotsc,n{-}r$.
 But this is just $|\pi|-|\pi|_{\widehat{\lambda}^t}$, which 
 proves~\eqref{Eq:critical}.
\end{proof}

\newtheorem{alemma}{Lemma}
\newtheorem{acor}[alemma]{Corollary}
\numberwithin{alemma}{section}

\appendix

\section{Root system miscellany}\label{A:root}

Our situation and notation will be as in the proof of 
Theorem~\ref{Theorem}.  To recap,
suppose that $G$ is a reductive algebraic group for which $G^{ss}$ is simple,
and let $P\subset G$ 
be a parabolic subgroup so that $G/P$ is cominuscule.
We will freely use the characterizations (i)---(iv) of cominuscule flag varieties
from Section~\ref{S:cominuscule}.
Fix a maximal torus $H$ of $P$ and let $L$ be the Levi subgroup of $P$ 
which contains $H$.  Let $v\in\frg/\frp$, and assume $v$ is neither
$0$, nor in the dense orbit of $L$ on $\frg/\frp$.
Let $Q\subset L$ be the stabilizer of $\frl\cdot v$, and define $\frz$ to 
be the quotient of the tangent space $\frg/\frp$ by $\frl\cdot v$.

We establish some essential facts about the root-space decomposition of the
Lie algebras $\frg, \frp, \frl, \frq,$ etc., as well as the 
subquotients $\frl\cdot v$ and $\frz$.  These results are needed 
in the proof of Theorem~\ref{Theorem}.
We begin with some general statements. 

Throughout, roots will mean the roots of $\frg$.
Let $\Phi$ be the set of roots of $\frg$, which are the non-zero
weights of $\frg$ under the action of the maximal torus $H$.
Once and for all, choose a non-zero
vector $v_\beta\in\frg_\beta$ in each
weight space of $\frg$.
If $\frs$ is a subquotient $H$-module of $\frg$, then we write
$\Phi(\frs)\subset\Phi$ for the non-zero weights of $\frs$.
If $\beta\in\Phi(\frs)$, then we also write 
$v_\beta\ (\in\frs)$ for the
image of $v_\beta\in\frg$ in $\frs$. 

If $\Delta\subset\Phi$ is a system of simple roots, then we may express any root
$\beta\in\Phi$ uniquely as an integral linear combination of the simple roots
in $\Delta$. 
Let $m_\delta(\beta)$ be the coefficient of $\delta \in \Delta$ in this 
expression for $\beta$.  
Write $m_\delta(\frg)$ for 
$\max_{\beta \in \Phi(\frg)} (m_\delta(\beta))$, which 
is the coefficient of $\delta$ in the highest root
of $\frg$. 

A root $\beta$ is \Blue{{\em positive}} (respectively \Blue{{\em negative}}) if any
coefficient $m_\alpha(\beta)$ is positive (respectively negative).
Since a root cannot be both positive and negative, we have the decomposition
$\Phi=\Phi^+\sqcup\Phi^-$ of $\Phi$ into positive and negative roots.
We say that a root $\beta$ is higher than $\beta'$ if $\beta-\beta'$
is a positive root. 
This definition depends upon the choice $\Delta$ of simple roots.
We say that $\Delta$ is \Blue{{\em compatible}} with $P$ if
$\Phi^-\subset\Phi(\frp)$. 
If $\alpha\in\Delta$ is the root defining $P$,
$\{\alpha\}=\Delta-\Phi(\frp)$,
then the weights of $\frl$ are exactly those $\beta\in\Phi$ such that
$m_\alpha(\beta)=0$. 

Recall that the standard pairing on the root space $\frh^*$ is
\[
   \langle \beta, \alpha\rangle\ :=\ 
    2\frac{(\beta,\alpha)}{(\alpha,\alpha)}\,.
\]
Here, $(\cdot,\cdot)$ is any $W$-invariant Euclidean inner product
on $\frh^*$.  By a long root, we mean any root $\alpha$ for which
$(\alpha,\alpha)$ is maximized.  If $G$ is simply laced, then every root
is long.

We recall the following basic facts about root systems.
Numbers 1 and 2 are found, for example, in Section 9.4 of~\cite{Hu72}.

\begin{enumerate}
\item[\bf 1.]
 If $\alpha$ is a long root then 
 $\langle \beta,\alpha\rangle\in\{-2,-1,0,1,2\}$
 and
 $\langle \beta,\alpha\rangle=\pm 2$ only if $\beta=\pm\alpha$.

\item[\bf 2.]
 If $\beta,\alpha\in\Phi$ with 
 $\pm\langle \beta,\alpha\rangle<0$, then $\beta\pm\alpha$ is a root.
 If $\alpha$ is a long root and $\beta+\alpha$ is a root, then 
 $\langle \beta,\alpha\rangle=-1$.

\item[\bf 3.]
 If a subgroup $K$ of $G$ contains the maximal torus $H$ and
 $\frs$ is a $K$-subrepresentation of $\frg$, then for every
 $\gamma\in\Phi(\frs)$ and $\beta\in\Phi(\frk)$
 with $\beta+\gamma\in\Phi$, we have
 $\beta+\gamma\in\Phi(\frs)$.  
 (Here, $\frk$ is the Lie algebra of $K$.)

 \noindent
 {\it Proof.} If $\beta+\gamma$ is a root then $v_\beta$ acts
 non-trivially on $v_\gamma$ and the result lies in
 $\frg_{\beta+\gamma}$, and so $\frg_{\beta+\gamma}\subset\frs$.
 \qed

\end{enumerate}

Given a system $\Delta\subset\Phi$ of simple roots of $\frg$,
a sequence
\[
  \gamma_1\ \to\ \gamma_2\ \to\ \dotsb\ \to\ \gamma_e
\]
of roots of $\frg$ is an \Blue{{\em increasing chain}} if, for all $k$,
$\gamma_{k+1}=\gamma_k+\delta_k$ where $\delta_k\in\Delta$.
That is, if at each step we raise by a simple root.
If $\gamma_1\in\Delta$, then for $\delta\in\Delta$, the coefficient
$m_\delta(\gamma_e)$ is the number of times $\delta$ was used in the chain
(including the first step $0 \to \gamma_1$).

\begin{alemma}\label{L:ONE}
 Let $K$ be any algebraic subgroup of $G$ containing $H$.
\begin{enumerate}
\item[(i)]
 If  $G/K$ is a cominuscule flag variety, then
 for any $\beta_1,\beta_2\in\Phi(\frg/\frk)$, 
 $\langle \beta_1,\beta_2 \rangle\geq 0$.
\item[(ii)]
 If $G/K$ is not a cominuscule flag variety, then there exist
 $\beta_1, \beta_2 \in \Phi(\frg/\frk)$ such that 
 $\beta_1 + \beta_2 \in \Phi \cup \{0\}$.
\end{enumerate}
\end{alemma}

Note that (i) implies the converse of (ii) and if 
$G$ is simply laced, then (ii) implies the converse 
of (i).
 
\begin{proof}
 Suppose that the homogeneous space $G/K$ is a cominuscule flag variety.
 Then in particular $K$ is a maximal parabolic subgroup.  
 Choose a system $\Delta$ of simple roots compatible with $K$, and let
 $\alpha\in\Delta$ be the simple root defining $K$. 
 Since $G/K$ is cominuscule, $m_\alpha(\beta)=1$ for every
 $\beta\in\Phi(\frg/\frk)$.
 Indeed, every root in $\Phi(\frg/\frk)$ lies in an increasing chain of roots 
 that starts with $\alpha$ and ends with the highest root.
 For (i), if $\langle\beta_1,\beta_2\rangle<0$ then
 $\beta_1+\beta_2\in\Phi(\frg/\frk)$ and so 
 $m_\alpha(\beta_1+\beta_2)=m_\alpha(\beta_1)+m_\alpha(\beta_2)=2$, 
 which is a contradiction.

 Suppose that $G/K$ is not cominuscule. 
 If $K$ is not a parabolic subgroup, then there exists a root $\gamma$
 of $\frg$ with neither
 $\gamma$ nor $-\gamma$ a root of $\frk$.  Thus we can take
 $\beta_1 = -\beta_2 = \gamma$.
 Otherwise, choose a positive system of roots 
 compatible with $K$, and let $\gamma_1$ be a simple root defining $K$. 
 Take an increasing chain of roots connecting $\gamma_1$ to the
 highest root,
\[
  \gamma_1\ \to\ \gamma_2\ \to\ \dotsb\ \to\ \gamma_{\rm top}\,.
\]
 Observe that each $\gamma_k\in\Phi(\frg/\frk)$.
 Since $G/K$ is not cominuscule, either there is another simple root in
 $\Phi(\frg/\frk)$ or else $m_{\gamma_1}(\gamma_{\rm top})\geq 2$.
 Thus at some point $\gamma_k$ in this chain, we will raise by a simple
 root $\delta\in\Phi(\frg/\frk)$.
 Thus $\gamma_{k+1}=\gamma_k+\delta$, so we can take $\beta_1 = \gamma_k$
 and $\beta_2 = \delta$.
\end{proof}

An \Blue{{\em orthogonal sequence of long roots}} in $\Phi(\frg/\frp)$ is
a sequence $\alpha_1, \ldots, \alpha_r \in \Phi(\frg/\frp)$, 
where $\alpha_i$ are long roots, and
$\langle \alpha_i, \alpha_j \rangle = 0 $ for $i \neq j$.
Such a sequence is maximal if every long root 
$\beta \in \Phi(\frg/\frp)$
is non-orthogonal to some $\alpha_i$. 
Orthogonal sequences of long roots play
a key role in the structure of $\frg/\frp$. 

If $G/P$ is cominuscule then every non-zero vector
$v \in \frg/\frp$ lies in the $L$-orbit of a sum
$v_{\alpha_1}+\dotsb+v_{\alpha_r}$, where
$\alpha_1,\dotsc,\alpha_r\in\Phi(\frg/\frp)$ is an orthogonal
 sequence of long roots~\cite{RRS92}.
Our assumption that $v$ does not lie in the dense orbit is equivalent
to assuming that 
$\alpha_1, \dotsc, \alpha_r$ is not maximal.
The construction of $Q$, $\frz$, etc. is $L$-equivariant with respect to
the choice of $v$, and thus we encounter no loss of generality in assuming
$v$ takes this normal form.  We therefore write
\[
  v\ =\ v_{\alpha_1}+\dotsb+v_{\alpha_r},
\] 
with
$\langle \alpha_i, \alpha_j \rangle = 0$ for $i\neq j$.  However, note
that in the following lemmas, whenever $\alpha_i$ does not appear 
explicitly in the statement, the result is valid for all non-zero 
$v \in \frg/\frp$ which are not in the dense orbit of
$L$.

\begin{alemma}\label{L:Claim1}
  If $\gamma\in\Phi(\frl)$, then there is at most one index $i$ such
  that $\langle \gamma,\alpha_i\rangle=-1$.
\end{alemma} 

\begin{proof}
 Suppose that $\langle \gamma,\alpha_i \rangle=-1$.
 Then $\gamma + \alpha_i \in \Phi$, and since $L$ 
 preserves $\frg/\frp$, we have 
 $\gamma + \alpha_i \in \Phi(\frg/\frp)$.
 If $j\neq i$, then Lemma~\ref{L:ONE} (ii) and
 $\langle \alpha_i,\alpha_j\rangle=0$ imply that
\[
   0\ <\ \langle\gamma+\alpha_i, \alpha_j\rangle \ =\ 
    \langle\gamma,\alpha_j\rangle\ .\ \  \qedhere
\]
\end{proof}

\begin{alemma}\label{L:rank1roots}
 Let $\alpha\in\Phi(\frg/\frp)$ be a long root.
 Then $\frl\cdot v_\alpha$ is $H$-invariant and 
 $\Phi(\frl\cdot v_\alpha)=\{\beta\in\Phi(\frg/\frp)\mid 
  \langle \beta,\alpha\rangle\geq 1\}$.
\end{alemma}

\begin{proof}
 As $G/P$ is cominuscule, all long roots of $\frg/\frp$ are conjugate (in
 fact by $W_L$~\cite{RRS92}), so we may
 assume that $\alpha\in\Delta$ defines $P$.
 Since $\frg_\alpha$  and $\frl$ are $H$-modules, 
 so is $\frl\cdot v_\alpha=\frl\cdot \frg_\alpha$.
  Note that  
\[
   \Phi(\frl\cdot v_\alpha)\ =\ 
     \{\beta\in\Phi(\frg/\frp)\mid \beta-\alpha\in \Phi(\frl)\cup\{0\}\}\,.
\]
 Let $\beta\in\Phi(\frg/\frp)$ and suppose that 
 $\langle \beta,\alpha\rangle\geq 1$.
 Recall that $m_\alpha(\beta)=1$.
 If  $\langle \beta,\alpha\rangle= 2$, then $\beta=\alpha$.
 Otherwise  $\langle \beta,\alpha\rangle= 1$ and so 
 $\beta-\alpha$ is a root.
 Then $m_\alpha(\beta-\alpha)=m_\alpha(\beta)-1=0$,
 and thus $\beta-\alpha\in\Phi(\frl)$. 

 We show the other inclusion.
 If $\beta-\alpha\in\Phi(\frl)\cup\{0\}$ then 
\[
  \langle \beta,\alpha\rangle\ =\ 
   \langle \beta-\alpha,\alpha\rangle+\langle \alpha,\alpha\rangle
   \ \geq\ -1 + 2\ =\ 1\,. \qedhere
\]
\end{proof}


Recall that $\frz:=(\frg/\frp)/(\frl\cdot v)$.

\begin{alemma}\label{L:frlcdotv}
 $\frl\cdot v$ is $H$-invariant and we have 
 \begin{eqnarray*}
  \Phi(\frl\cdot v)&=& \{\beta\in\Phi(\frg/\frp)\mid 
            \langle \beta,\alpha_i\rangle\geq 1\,,
     \text{ for some $i$}\}, \ \text{and}\\
  \Phi(\frz)&=& \{\beta\in\Phi(\frg/\frp)\mid 
            \langle \beta,\alpha_i\rangle=0\,,
     \text{ for all $i$}\}\,.
 \end{eqnarray*}
%
%
This holds even when $v$ lies in the dense orbit of $L$.
\end{alemma}

\begin{proof}
 We claim that
 \[
   \frl\cdot v\ =\ \frl\cdot v_{\alpha_1}+\frl\cdot v_{\alpha_2}+\dotsb +
   \frl\cdot v_{\alpha_r}\,,
 \] 
  from which the statement of the lemma follows from
  Lemma~\ref{L:rank1roots}, as each 
 $\alpha_i$ 
 is a long root.
  For each $i=1,\dots,r$, let $\frl_i\subset\frl$ be the linear span of the
  set
\[
   \Gamma_i\ :=\ \{ v_\gamma\mid \gamma\in\Phi(\frl)\mbox{\ and\ }
          \langle\gamma,\alpha_i\rangle=-1\}\,.
\]
 Then we have 
 $\frl\cdot v_{\alpha_i}=\frl_i\cdot v_{\alpha_i}=\frl_i\cdot v$. 
 The last equality is a consequence of Lemma~\ref{L:Claim1}.
 Lemma~\ref{L:Claim1} also implies that the sets $\Gamma_i$ are disjoint and
 therefore $\frl_1+\frl_2+\dotsb+\frl_r$ is a direct sum.
 Thus we have
\begin{eqnarray*}
  \frl\cdot v &=& \frl\cdot(v_{\alpha_1}+v_{\alpha_2}+\dotsb+v_{\alpha_r})\\
              &\subset& \frl\cdot v_{\alpha_1}+\frl\cdot v_{\alpha_2}
                       +\dotsb+\frl\cdot v_{\alpha_r}\\
     &=& \frl_1\cdot v_{\alpha_1}+\frl_2\cdot v_{\alpha_2}
                       +\dotsb+\frl_r\cdot v_{\alpha_r}\\
     &=& \frl_1\cdot v+\frl_2\cdot v+\dotsb+\frl_r\cdot v\\
     &=& (\frl_1+\frl_2+\dotsb+\frl_r)\cdot v\ \subset\ \frl\cdot v\,,
\end{eqnarray*}
which proves the claim.
\end{proof}

\begin{alemma}\label{L:2orthroots}  
 Let $\alpha_1, \alpha_2, \dotsc \alpha_r \in \Phi(\frg/\frp)$ be an
 orthogonal sequence of long roots.
\begin{enumerate}
\item[(i)]
For any $\beta \in \Phi(\frg/\frp)$, 
 there exist at most two distinct indices $i$ such that 
 $\langle\beta,\alpha_i\rangle\ \geq 1$.

\item[(ii)] There exists $\beta \in \Phi(\frg/\frp)$ such that 
$\langle \beta, \alpha_1 \rangle 
= \langle \beta, \alpha_2 \rangle = 1$, when $r\geq 2$.
\end{enumerate}
\end{alemma}

\begin{proof}
 Choose a system $\Delta$ of simple roots compatible with $P$ and let
 $\alpha\in\Delta$ be the simple root defining $P$. 

 (i) Suppose there are three indices, say $i, j, k$.  Then
 $\langle\beta,\alpha_{i}\rangle\geq 1$ so $\beta-\alpha_{i}$
 is a root.  
 Then $\langle\beta-\alpha_{i},\alpha_{j}\rangle\ \geq 1$, so
 $\beta-\alpha_{i}-\alpha_{j}$ is a root.  Similarly 
 $\beta-\alpha_{i}-\alpha_{j}-\alpha_{k}$ is a root.  But
 now 
 $m_\alpha(\beta-\alpha_{i}-\alpha_{j}-\alpha_{k})=-2$ and
 there is no root with this property as $G/P$ is cominuscule.

 (ii) Since $W_L$ acts transitively on all orthogonal sequences of 
long roots of the same length~\cite{RRS92}, it suffices to show this for
a particular pair of orthogonal long roots.  
Set $\alpha_1 := \alpha$, the simple root defining $P$, 
and let $\alpha_2$ be the highest root of $\frg$.
If an orthogonal pair of long roots exists, Lemma~\ref{L:Simple} implies 
this is such a pair (and the argument is non-circular).
Let $\delta\in\Delta$ be a root such that
$\langle\delta,\alpha_2\rangle=1$
and consider the sum, $\beta$, of $\alpha+\delta$ with all the simple
roots in the Dynkin diagram of $G$ which lie strictly between $\alpha$
and $\delta$.  
Such a sum is always a root.
We have 
$\langle \beta, \alpha_1 \rangle 
= \langle \beta, \alpha_2 \rangle = 1$.
\end{proof}

Recall that $Q=\Stab_L(\frl\cdot v)$ so that $\frq$ is spanned by those
$v_\gamma$ which stabilize $\frl\cdot v$.

\begin{alemma}\label{L:DeltaQ}
 $\Phi(\frq)=\{\gamma\in\Phi(\frl)\mid
   \langle\gamma,\alpha_i\rangle= 0\text{ for all $i$ \Blue{\em or} }
   \langle\gamma,\alpha_i\rangle\geq 1\text{ for some $i$}
  \}$.
\end{alemma}

 Observe that by Lemma~\ref{L:Claim1} we deduce,
\[
   \Phi(\frl/\frq)\ =\ 
   \{\gamma\in\Phi(\frl)\mid
   \langle\gamma,\alpha_i\rangle\leq 0\text{ for all $i$ \Blue{{\em and}} }
   \langle\gamma,\alpha_i\rangle =-1\text{ for exactly \Blue{{\em one}} $i$}
  \}\,.
\]

\begin{proof}
 Let $\gamma\in\Phi(\frl)$ and suppose that 
 $\langle\gamma,\alpha_i\rangle\geq 0$, for all $i=1,\dotsc,r$.
 If $\beta\in\Phi(\frl\cdot v)$, then Lemma~\ref{L:frlcdotv} implies
 that $\langle \beta,\alpha_i\rangle\geq 1$ for some $i$,
 and therefore $\langle\gamma+\beta,\alpha_i\rangle \geq 1$.
 Similarly, suppose that 
 $\langle\gamma,\alpha_i\rangle\geq 1$ for some index $i$.
 If $\beta\in\Phi(\frl\cdot v)\,(\subset \Phi(\frg/\frp))$, 
 then Lemma~\ref{L:ONE}(i) implies that 
 $\langle\beta,\alpha_i\rangle \geq 0$, and so again we have 
 $\langle\gamma+\beta,\alpha_i\rangle \geq 1$.
 Thus in either case, Lemma~\ref{L:frlcdotv} implies that if $\gamma+\beta$
 is a root, then it lies in $\Phi(\frl\cdot v)$, and so $v_\gamma\in\frq$
 as it stabilizes $\frl\cdot v$.

 For the converse, suppose
 that $\langle\gamma,\alpha_i\rangle \leq 0$ for all
 $i$, and $\langle\gamma,\alpha_j\rangle =-1$ for some index $j$ with
 $1\leq j\leq r$.
 Suppose moreover that $v_\gamma$ stabilizes $\frl\cdot v$.
 We show this leads to a contradiction.\smallskip

\noindent
 {\bf Claim.} Assume that $\beta\in\Phi(\frg/\frp)$ is not equal to $\alpha_j$.
 Then $\langle\beta,\gamma\rangle<0$ only if there exist exactly
 two indices $i$ such that $\langle\beta,\alpha_i\rangle \geq 1$.\smallskip

\noindent
 {\em Proof of Claim.}  
 Suppose that $\langle \beta,\gamma\rangle<0$.
 Then $\gamma+\beta$ is a root, and since $\frg/\frp$ is an $\frl$-module, 
 $\gamma+\beta\in\Phi(\frg/\frp)$.
 By Lemma~\ref{L:ONE}(i) we have 
 $0 \leq \langle\gamma+\beta,\alpha_j\rangle\ = 
 -1+ \langle\beta,\alpha_j\rangle$.
 Since $\beta \neq \alpha_j$, we have 
 $\langle \beta,\alpha_j\rangle=1$, 
 and so $\beta\in\Phi(\frl\cdot v)$, by Lemma~\ref{L:frlcdotv}.

 Since $v_\gamma$ stabilizes $\frl\cdot v$, we must have
 that $\gamma+\beta\in \frl\cdot v$, and thus there is some index $i$ with
 $\langle \gamma+\beta,\alpha_i\rangle\geq 1$.
  Then
\[
  1\ \leq\ \langle \gamma+\beta,\,\alpha_i\rangle\ =\ 
  \langle \gamma,\,\alpha_i\rangle\ +\
  \langle \beta,\,\alpha_i\rangle \ \leq\ 
  \langle \beta,\,\alpha_i\rangle\,,
\]
 so $1\leq\langle\beta,\alpha_i\rangle$.
 Necessarily, $i\neq j$ as  $\langle\gamma+\beta,\alpha_j\rangle
 = \langle\gamma,\alpha_j\rangle + \langle\beta,\alpha_j\rangle 
 =0$.
\qed\smallskip

 Now since $L\cdot v$ is not dense, there exists a long root 
 $\alpha \in \Phi(\frg/\frp)$ orthogonal to $\alpha_i$ 
 for all indices $i$.
 We may assume that $\alpha$ is the simple
 root defining $P$.

 Consider the set $h := \{\beta \in \Phi(\frg/\frp) \mid
 \langle \beta, \alpha\rangle > 0 \}$.  By Lemma~\ref{L:2orthroots}(i) 
 if $\beta \in h$ then we have $\langle \beta, \alpha_i\rangle>0$
 for at most one index $i$ ($\beta$ is already positively paired
 with $\alpha$).
 Also, $\alpha_j\not\in h$, as 
 $\langle\alpha_j,\alpha\rangle=0$.
 Thus by the claim, we have 
 $\langle \beta, \gamma \rangle \geq 0$ for all $\beta \in h$.

 Now consider the sum, $\beta_0$, of $\alpha$ and all the simple roots in
 the Dynkin diagram of $G$ which lie strictly between $\alpha$ and the
 nearest simple root used in $\gamma$.
 Such a sum is always a root, and $\beta_0\in\Phi(\frg/\frp)$.
 Moreover, $\langle\beta_0, \alpha\rangle =1$,
 so $\beta_0 \in h$, but $\langle\beta_0, \gamma\rangle<0$, a
 contradiction.
\end{proof}

In summary, the roots of $\frg$ decompose into a disjoint union of
\[
  \Phi(\frr) = \Phi(\frn_P) \sqcup \Phi(\frq),\ \  \Phi(\frl/\frq),\ \  
   \Phi(\frl\cdot v),\ \ \mbox{\rm and}\ \ \Phi(\frz)\,.
\]
The roots, $\gamma$, of these pieces are characterized by their pairings with
respect to the long roots $\alpha_i$ and the values of $m_\alpha(\gamma)$,
where $\alpha$ is the simple root defining $P$.
These characterizations are given concisely in Table~$\ref{Table:decompsummary}$.

\begin{table}[htbp]
\centering
\begin{tabular}{|c|c|l|}
\hline
$\gamma $ & $m_\alpha(\gamma)$ & 
 $\langle \gamma, \alpha_i\rangle$ \\
\hline
\hline
$\Phi(\frz)$ & $1$ & $=0$ for all $i$ 
   \\ \hline
$\Phi(\frl\cdot v)$ & $1$ & $\geq 0$ for all $i$\ \ and\ \ $\geq 1$ for some $i$ 
   \\ \hline
$\Phi(\frl/\frq)$ & $0$ & $\leq 0$ for all $i$\ \ and\ \ $=-1$ for exactly one $i$ 
   \\ \hline
$\Phi(\frq)$ & $0$ & $= 0$ for all $i$\ \ or\ \ $\geq 1$ for some $i$ 
   \\ \hline
$\Phi(\frn_P)$ & $-1$ & $\leq 0$ for all $i$ 
   \\ \hline
\end{tabular}\vspace{4pt}
\caption{Summary of decomposition of $\Phi$}\vspace{-18pt}
\label{Table:decompsummary} 
\end{table}

\begin{alemma}\label{L:Qcominuscule}
 The subgroup $Q$ is a parabolic subgroup of $L$, and the flag 
 variety $L/Q$ is cominuscule.
\end{alemma}

\begin{proof}
 Suppose that $L/Q$ is not cominuscule.
 By Lemma~\ref{L:ONE}(ii), there must be
 be two roots $\gamma_1,\gamma_2\in\Phi(\frl/\frq)$ with 
 $\gamma_1+\gamma_2 \in \Phi(\frl) \cup \{0\}$
 
 By Lemma~\ref{L:DeltaQ}, 
 $\langle\gamma_1,\alpha_i\rangle, \langle\gamma_2,\alpha_i\rangle\leq 0$
 for all $i$ and there exist unique indices $i_1$ and $i_2$ such that
\[
  \langle \gamma_1,\alpha_{i_1}\rangle \ =\  -1
   \text{ \ and \ }
  \langle \gamma_2,\alpha_{i_2}\rangle \ =\  -1\,.
\]
 If $\gamma_1 = -\gamma_2$,
 then $\langle\gamma_2,\alpha_{i_1}\rangle = 1$, which contradicts
 Lemma~\ref{L:DeltaQ}.
 Otherwise,
 $\gamma_1+\gamma_2$ is a root, necessarily in $\Phi(\frl/\frq)$.
 If $i_1=i_2$ then  
 $\langle \gamma_1+\gamma_2,\alpha_{i_1}\rangle =-2$, which is
 impossible.
 Otherwise we have 
 $\langle \gamma_1+\gamma_2,\alpha_{i_1}\rangle \leq-1$ and
 $\langle \gamma_1+\gamma_2,\alpha_{i_2}\rangle \leq-1$,
 which is also impossible, by Lemma~\ref{L:Claim1}.
\end{proof}

Let $R$ be the parabolic subgroup of $P$ corresponding to $Q \subset L$
so that $R=\Stab_P(\frl\cdot v)$.
Since $\frl\cdot v$ is $H$-stable, $H\subset Q\subset R$.
We also assume that our system $\Delta$ of simple roots is
compatible with the parabolic subgroups $P$ and $R$. 

\begin{alemma}\label{L:Simple}
  The highest root of $\frg$ is an element of $\Phi(\frz)$.
  If $\gamma\in\Phi(\frl)$ is a simple root defining $Q$, then 
  $m_\gamma(\beta) = m_\gamma(\frg)$ for all
  $\beta\in\Phi(\frz)$, and $\frz$ is an irreducible $L_Q$-module.
  The same statements hold for $R$ in place of $Q$.
\end{alemma}

\begin{proof}
 Let $\beta_1\in\Phi(\frz)\subset\Phi(\frg/\frp)$ and consider an
 increasing chain of roots
 \begin{equation}\label{E:z-chain}
   \beta_1\ \to\ \beta_2\ \to \ \dotsb\,.
 \end{equation}
 Let $\delta_k$ be the simple root $\delta_k:=\beta_{k+1}-\beta_k$.
 Let $\alpha$ be the simple root defining $P$.
 Then
 $1 = m_\alpha(\beta_1) \leq m_\alpha(\beta_k) \leq m_\alpha(\frg) = 1$.
 Hence $\beta_k \in \Phi(\frg/\frp)$ 
 and $\delta_k \neq \alpha$,
 and we conclude that $\delta_k$ is a simple root of $L$.  

Suppose $\delta_k$ is a simple root of $L_Q$.  Since $\frl \cdot v$
is a $Q$-submodule of $\frg/\frp$, we have the decomposition
$\frg/\frp = (\frl \cdot v) \oplus \frz$ as $L_Q$-modules (since
$L_Q$ is reductive).  In particular $\frz$ is an $L_Q$-module.
Thus if $\beta_k \in \Phi(\frz)$, then
$\beta_{k+1} = \beta_k + \delta_k \in \Phi(\frz)$.

 Otherwise, $\delta_k=\gamma$, a simple root of $L$ defining $Q$.
 Then $m_\gamma(\beta_{k+1})=m_\gamma(\beta_k)+1$.
 By Lemma~\ref{L:DeltaQ}, there is some index $i$ such that 
 $\langle\gamma,\alpha_i\rangle=-1$, as $\gamma\in\Phi(\frl/\frq)$.
 Then Lemma~\ref{L:frlcdotv} implies that 
 $\langle \beta_{k+1}, \alpha_i\rangle=
  \langle \beta_{k},   \alpha_i\rangle-1=-1$,
 which contradicts Lemma~\ref{L:ONE}(i).

 We conclude that every simple root $\delta_k$ arising from our
 chain~\eqref{E:z-chain} is a simple root of $Q$, and each
 $\beta_k\in\Phi(\frz)$. 
 This implies that $m_\gamma(\beta_k)$ is a constant, where $\gamma$ is
 a root of $L$ defining $Q$.
 Since every root of $\Phi(\frz)$ may be connected to the highest root
 of $\Phi(\frg/\frp)$, this highest root lies in $\Phi(\frz)$, 
 $\frz$ is an irreducible representation of $L_Q$, and 
 $m_\gamma(\beta)$ is constant for  $\beta\in\Phi(\frz)$, where $\gamma$
 is a root of $L$ defining $Q$, and this constant value is
 $m_\gamma(\frg)$.
\end{proof}

\begin{acor}\label{C:frz}
 We have the $R$-module isomorphism $Z(\frn_R)^*\simeq \frz$.
\end{acor}

\begin{proof}
 The dual $Z(\frn_R)^*$ of the center of $\frn_R$ is spanned by the
 vectors $v_\beta$ for which $m_\gamma(\beta) = m_\gamma(\frg)$
 for all simple roots $\gamma\in\Phi(\frn_R)$.
 Thus by Lemma
 \ref{L:Simple} we have an injective $R$-module morphism from 
 $\frz$ to $Z(\frn_R)^*$.  Since both $\frz$ and $Z(\frn_R)$ are 
 irreducible $L_R$-modules, this is an isomorphism.
\end{proof}

Let $G'$ be the subgroup $G'\ := \ Z_G ( Z_H ( Z(N_R)))$,
$P':= G'\cap R$, and let $L'$ be the Levi subgroup of $P'$.
Note that $G'\supset H$, so that $G'$ is determined by the weights
$\Phi(\frg')$ of its Lie algebra.

\begin{alemma}\label{L:roots}
 $\Phi(\frg')= \mathbb{Q}\Phi(\frz)\cap \Phi(\frg) \subset
\{\gamma \in \Phi(\frg)\ |\ \langle \gamma, \alpha_i \rangle
   = 0 \text{ for $i=1,\dotsc, r$}\}$.
\end{alemma}

\begin{proof}
First note that 
$Z_\frh(Z(\frn_R)) = \Phi(Z(\frn_R))^\perp = \Phi(\frz)^\perp$.  
Also, for any subalgebra $\frh' \subset \frh$,
$\Phi(Z_\frg(\frh')) = (\frh')^\perp \cap \Phi(\frg)$.
Thus
 \[
   \Phi(\frg')\ =\  ( \Phi(\frz)^\perp)^\perp\ \cap \ \Phi(\frg)
  \ =\  \mathbb{Q}\Phi(\frz)\cap \Phi(\frg)\,,
\]
 proving the equality.
 The inclusion is a consequence of Lemma~\ref{L:frlcdotv}.
\end{proof}

\begin{alemma}\label{L:stabCv}
$L' \subset Z_H(\frz) \Stab_L(\CC v)$.
\end{alemma}

\begin{proof}
First note that $\Phi(\stab_\frl(\CC v))
 = \{\gamma \in \Phi(\frl)\ |\ \langle \gamma, \alpha_i \rangle
 = 0 \text{ for $i=1,\dotsc, r$}\}$, which contains 
 $\Phi(\frg' \cap \frl) = \Phi(\frl')$, by Lemma~\ref{L:roots}.  
 Thus it suffices to show that 
$H \subset Z_H(\frz) \Stab_H(\CC v)$.  But since $\QQ\phi(\frz)$ and
$\QQ\{\alpha_1, \dotsc, \alpha_r\}$ are orthogonal, their
annihilators,
$\Phi(\frz)^\perp = Z_\frh(\frz)$ and $\{\alpha_1, \dotsc, \alpha_r\}^\perp
\subset \stab_\frh(\CC v)$ together span $\frh$.
\end{proof}

\begin{alemma}\label{L:sspart}
 The semisimple part of $G'$ is simple.
\end{alemma}

\begin{proof}
For any subset $\Gamma \subset \Phi(\frg')$, we form a graph by joining
$\gamma_1, \gamma_2 \in \Gamma$ by an edge if 
$\langle \gamma_1, \gamma_2 \rangle \neq 0$.
It suffices to show that there is a subset $\Gamma$ of $\Phi(\frg')$
which spans $\QQ\Phi(\frg')$ such that this graph is connected.
We show that $\Gamma = \Phi(\frz)$ is such a subset.

Extend
$\alpha_1, \dotsc, \alpha_r$ to a maximal orthogonal sequence
of long roots $\alpha_1, \dotsc, \alpha_m$.
By Lemma~\ref{L:2orthroots}(ii)
any pair $\alpha_i, \alpha_j \in \Phi(\frz)$ have a common
non-orthogonal root $\beta\in\Phi(\frg/\frp)$. 
By Lemma~\ref{L:2orthroots}(i), 
$\beta$ is orthogonal to $\alpha_1, \dotsc, \alpha_r$, hence
in $\Phi(\frz)$. 
Finally, every root $\beta \in \Phi(\frz)$ is 
non-orthogonal to some $\alpha_i$ (necessarily in $\Phi(\frz)$).
Indeed, as $\alpha_1,\dotsc,\alpha_m$ is maximal,
$\frg/\frp=\frl\cdot(v_{\alpha_1}+\dotsb+v_{\alpha_m})$.
Then, by Lemma~\ref{L:frlcdotv},
\[
   \{\beta \in \Phi(\frg/\frp)\ |\ \langle \beta, \alpha_i \rangle \geq  1
   \text{ for some $i$}\} 
   \ =\ \Phi(\frl \cdot (v_{\alpha_1} + \dotsb + v_{\alpha_m})) 
   \ =\ \Phi(\frg/\frp)\,.\qedhere
\]
\end{proof}


\begin{alemma}\label{L:nilr}
  The nilradical $\frn_{P'}$ is equal to $Z(\frn_R)$, the center of the
  nilradical of $R$.  In particular $G'/P'$ is cominuscule.  
\end{alemma}

As the Killing form on $\frg'$ identifies $(\frn_{P'})^*$ with $\frg'/\frp'$, 
this identifies $\frz$ with $\frg'/\frp'$.

\begin{proof}
 By our definition of $P'$, $\frn_{P'}\subset\frn_R$, and 
 $\frn_{P'}$ is an $H$-module.
 By Lemma~\ref{L:roots}, $Z(\frn_R)=\frz^*\subset \frn_{P'}$.
 Let $\gamma\in\Phi(\frn_{P'})$ be a weight that is not a weight of
 $Z(\frn_R)$.
 Then $-\gamma$ is either in $\Phi(\frl/\frq)$ or else in $\Phi(\frl\cdot v)$,
 and thus there is some $i=1,\dotsc,r$ such that 
 $\langle -\gamma,\alpha_i\rangle\neq 0$, by Table 2.
 In particular, $-\gamma\not\in\QQ\Phi(\frz)$, and so is not a 
 weight of $G'$.  As 
 $N_P' = Z(N_R)$ is abelian, we deduce that $G'/P'$ is cominuscule.
\end{proof}

\begin{alemma}\label{L:aseffective}
For every $q \in R$, there exists 
$l \in L' \cap \Stab_L(\CC v)$ such that 
for every $z \in \frz$ we have $qz = lz$.
\end{alemma}

\begin{proof}
 First, we show that if
 $\gamma\in\Phi(\frr)-\Phi(\frp')$, then the weight
 vector $v_\gamma$ acts trivially on $\frz$.
 This weight $\gamma$ does not lie in $\QQ\Phi(\frz)$, and thus
 if $\beta\in\Phi(\frz)$, then 
 $\beta+\gamma\not\in\Phi(\frz)$.
 Since $\frz$ is an $\frr$-module, this implies that 
 $v_\gamma\cdot v_\beta=0$.

It follows that there exists $p \in P'$ such that $pz = qz$
for every $z \in \frz$.
However, as $G'/P'$ is cominuscule, $N_{P'}$ acts trivially on
its Lie algebra, and hence acts trivially on $\frz$, so we can 
replace $p$ by an element of $L'$.  Finally 
as $Z_H(\frz)$ acts trivially on $\frz$, by Lemma~\ref{L:stabCv}
we can reduce further to $L' \cap \Stab_L(\CC v)$.
\end{proof}

Let $v' \in \frg'/\frp'=\frz$.  We assume that $v'$ is of the form
$v' = v_{\alpha_{r+1}} + \cdots + v_{\alpha_{r+r'}}$, where
$\alpha_{r+1}, \ldots, \alpha_{r+r'}$ is an orthogonal sequence of long
roots.
Set $v_1 := v+v'$.
Define $Q' := \Stab_{L'}(\frl \cdot v')$, and 
$Q_1 = \Stab_L(\frl \cdot v_1)$, and let $R'$ and $R_1$ be the
corresponding parabolic subgroups in $P'$ and $P$, respectively.

\begin{alemma}\label{L:Q1standard}
$Q' = L' \cap Q_1$.
\end{alemma}

\begin{proof}
 This follows from the characterization of $\Phi(\frq')$ and
 $\Phi(\frq_1)$ of Lemma~\ref{L:DeltaQ} as
 $\langle \gamma_1, \alpha_i \rangle = 0$, for $\gamma\in\Phi(\frl')$ and
 each $i=1, \dotsc, r$.
\end{proof}

\begin{alemma}\label{L:indfrn}
 $Z(\frn_{R_1})=Z(\frn_{R'})$.
\end{alemma}

\begin{proof}
Note that the weights of $Z(\frn_{R_1})$ are exactly those in 
$\Phi(\frn_{P})$ which are annihilated by $\alpha_{i}$ for $i=1,\dotsc,r'$.
The weights in $Z(\frn_{R'})$ are weights in $\Phi(\frn_{P'})$ which are 
annihilated 
$\alpha_{i}$ for $i=r+1,\dotsc,r'$.  
Since $\Phi(\frn_{P'})$ are the weights of $\frn_P$ annihilated by
$\alpha_{i}$ for $i=1,\dotsc,r$, the result follows.
\end{proof}


\providecommand{\bysame}{\leavevmode\hbox to3em{\hrulefill}\thinspace}
\providecommand{\MR}{\relax\ifhmode\unskip\space\fi MR }
\providecommand{\MRhref}[2]{%
  \href{http://www.ams.org/mathscinet-getitem?mr=#1}{#2}
}
\providecommand{\href}[2]{#2}

\end{document}